\pdfoutput=1
\documentclass[11pt,a4paper]{article}
\usepackage{amsmath,amssymb,amsthm}
\usepackage{geometry}
\usepackage{bm}
\usepackage{graphicx}
\usepackage[T1]{fontenc}
\usepackage{inconsolata}  
\usepackage{booktabs}
\usepackage{xcolor}
\usepackage{placeins}   

\usepackage[numbers]{natbib}
  \usepackage{doi}        
  \usepackage{hyperref}

\newcommand{\numm}[1]{\texttt{#1}}
\newcommand{\numb}[1]{\textbf{\texttt{#1}}}

\newtheorem{remark}{Remark}

\newcommand{\kT}{k_{\mathrm{B}}T}
\newcommand{\gU}{\nabla U}
\newcommand{\hf}{\hat{f}}
\newcommand{\R}{\mathbb{R}}
\providecommand{\varpsi}{\psi}

\title{Accelerated sampling using SamAdams variable timesteps
       and position-adaptive Langevin dynamics}
\author{Benedict Leimkuhler\thanks{School of Mathematics, University of Edinburgh, Edinburgh, UK.} 
\and Peter A. Whalley\thanks{Seminar for Statistics, Department of Mathematics, ETH Z\"urich, Z\"urich, Switzerland.}}
\date{}

\begin{document}
\maketitle

\begin{abstract}
We introduce an accelerated Langevin-based sampling method that is based on two complementary devices: \emph{SamAdams} adaptive timestepping, which automatically shrinks the effective integration step in stiff regions of phase space using a relaxed stiffness monitor, and \emph{position-adaptive Langevin} (PAL) dynamics, which concentrates friction along the local force direction while preserving the canonical distribution as the exact invariant measure. The resulting combined scheme (SA-PAL) is implemented in a palindromic integrator which requires only one force evaluation per iteration through suitable organisation of the integration steps and by exploiting the rank-one-plus-scalar structure of the PAL friction tensor.  We test the method on various model problems: the Rosenbrock function, a thin entropic channel,  the Mueller-Brown potential, and a Bayesian parameterisation problem with a sparsity-inducing shrinkage prior.  On the Rosenbrock and Mueller-Brown potentials mixing rates are improved by 1.5-3 times compared to fixed stepsize integration. Efficiency gains of more than an order of magnitude are documented in the other examples.  
\end{abstract}
\section{Introduction}
\label{sec:intro}

Underdamped Langevin dynamics on $\mathbb{R}^{2d}$ with constant symmetric, positive-definite mass matrix $M \in \mathbb{R}^{d\times d}$ and friction matrix $\Gamma \in \mathbb{R}^{d\times d}$, is defined by the following SDE system:
\begin{align}
  dx &= M^{-1}p\,dt, \label{eq:UL-x}\\
  dp &= -\gU(x)\,dt - \Gamma M^{-1} p\,dt + \sqrt{2\kT\Gamma}\,dW_t, \label{eq:UL-p}
\end{align}
where $U:\mathbb{R}^{d} \to \mathbb{R}$ is a potential, $(W_t)_{t\geq 0}$ denotes a vector of independent Wiener processes on $\mathbb{R}^{d}$, $T$ is the temperature and $k_{B}$ is Boltzmann's constant.
Under mild assumptions, the invariant measure of \eqref{eq:UL-x}--\eqref{eq:UL-p} is the canonical distribution with density defined pointwise up to a normalisation constant by
$\pi(x,p) \propto \exp(-H(x,p)/\kT)$, where the Hamiltonian  is defined via $H(x,p) = U(x)+\tfrac12 p^T M^{-1} p$,
for any constant positive-definite $\Gamma$. Additionally, under mild assumptions, the system is ergodic with respect to the invariant measure \cite{pavliotis2014}. As a consequence, it is a popular method for sampling in molecular modelling and Bayesian computation. 

The practical efficiency of \eqref{eq:UL-x}-\eqref{eq:UL-p} depends crucially on the choice of friction: if the friction is too weak, trajectories exhibit long ballistic excursions and may explore inefficiently, whereas if friction is too strong, momentum is rapidly overdamped and the dynamics approaches the slower overdamped Langevin regime \cite{cao2023explicit,lu2026sharp}.  The familiar scalar choice $\Gamma=\gamma I$ damps every direction equally.  Introducing a friction matrix and allowing it to vary with position, $\Gamma: \R^d \rightarrow \R^{d\times d}$, instead lets different directions be damped independently, and in a way that may be tuned to local geometric features of the landscape under study.

Many potentials of interest involve both very stiff directions (large curvature) and flat or slowly varying directions.  In case a fixed stepsize is used, it must be small enough to resolve the stiffest
mode \cite{LePaWh2024}, sacrificing efficiency in the flat regions where much larger steps would
be safe.  Adaptive (variable) timestepping integrators address this by contracting
the step only where necessary, but doing so in a way that preserves the
invariant measure is non-trivial. 

In this work we propose an algorithm that combines variable timesteps with a specific position-dependent friction tensor as follows:
\begin{enumerate}
\item\textbf{SamAdams variable timestep~\cite{samadams}.}  The integrator continuously
  senses how demanding the local dynamics are through the magnitude of the
  force, or in tighter geometries through the speed of the motion, and smoothly
  contracts the integration step wherever the demand is high, relaxing it again
  where possible.  Because the effective step never exceeds a chosen nominal
  value, that nominal step can be set well beyond the stability threshold of a
  conventional fixed-stepsize integrator.  The accompanying
  time-rescaling is corrected for in a simple reweighting procedure, so that Monte Carlo averages from the canonical distribution can be recovered.  We use two variants of the local sensor, one driven by the
  force magnitude (the recommended default), and one that also responds to how fast the trajectory is changing, better suited to a narrow entropic bottleneck.  Both options are defined in Section~\ref{sec:samadams}.
\item \textbf{Position-adaptive Langevin (PAL).}  The friction tensor is taken to be the normalized rank-one-plus-scalar form
  \begin{equation}
    \Gamma(x) = \alpha I + \beta\,\hf(x)\hf(x)^T,
    \qquad \hf(x) = \frac{\gU(x)}{|\gU(x)|},
    \label{eq:PAL-friction}
  \end{equation}
  with constants $\alpha > 0$ and $\beta > -\alpha$.  The eigenvalues of $\Gamma$ are the spatially constant values $\alpha$ (transverse to $\gU$) and $\alpha+\beta$ (parallel to $\gU$); only the eigendirection rotates with $x$.  This structure admits an $O(d)$ exact solution of the Ornstein--Uhlenbeck (OU) step and guarantees that the canonical distribution is the exact invariant measure.
\end{enumerate}
   
 The idea of using a variable friction is not entirely new even if it is rarely used in practice.  The paper of Lim and Tao ~\cite{LimTao2025} considered this and the specific choice of the Hessian of the potential as the friction matrix and showed theoretically that the continuous dynamics could enhance performance (even if no numerical experiments were provided in that paper).

Due to the special form (scaled identity + rank-1 for friction, SamAdams variable stepsizes), the scheme we describe is also practically free of computational overhead compared to standard integrators.   The two devices are combined in the AZBOBZA palindromic integrator described herein (Section~\ref{app:algorithm}), which requires a single force evaluation per step, the same cost as other typical 2nd order Langevin integrators \cite{LeMa2013}, and passes all computed quantities through each sub-step without redundant work other than a few extra inner products.

In SA-PAL, the two devices address the temporal and directional components simultaneously, and this is shown in the numerical experiments of Section~\ref{sec:numerics}.  It is observed that the PAL contribution is most visible in landscapes with clearly curved low-energy structure, such as the Rosenbrock potential. On the thin channel, much of the gain comes from the SamAdams timestep; here the fixed-step chains barely cross the corridor at all, so SA-PAL's main effect is to pass through it far more often, although in experiments the position adaptive friction also provides a quantifiable benefit.

On the relatively smooth Mueller-Brown potential, the picture is more intricate and reveals a further benefit of the new scheme. Compared against the fastest-mixing isotropic friction the gain in position exploration looks modest, but at that friction, it turns out that the fixed stepsize solution is substantially biased: it mixes the positions quickly but crosses between the three wells too rarely to sample their relative populations correctly. The proper comparison is against the isotropic friction that matches SA-PAL's accuracy, and in relation to this choice SA-PAL mixes the positions about $2.8\times$ faster while reducing the potential-energy autocorrelation several times more quickly. The benefit PAL adds on top of SA is a directional decoupling of damping rates; it damps the force direction which drives the variation of energy without over-damping the slow inter-well motion. Thus, in this example, SA-PAL achieves a more favourable balance between exploration speed and equilibrium accuracy than any of the isotropic-friction baselines considered. This improvement in both mixing and accuracy is visible whenever the landscape has sufficient complexity.

Beyond these low-dimensional examples, we also apply SA-PAL to a realistic, higher-dimensional target, namely, Bayesian linear regression under a horseshoe prior, using an example borrowed from \cite{Carvalho2010}, whose posterior has the sharp ``funnel'' geometry that is notoriously difficult for fixed-step samplers. Here the combined scheme delivers order-of-magnitude reductions in integrated autocorrelation time relative to a well-tuned fixed-step baseline while keeping the estimated posterior moments accurate. 

\paragraph{Related work.}
Adaptive-friction Langevin samplers in the literature can be grouped into several distinct families, and the rank-one-plus-scalar friction \eqref{eq:PAL-friction} adopted here does not coincide with any of them.

Chak et al \cite{Chak2023} analyze the optimal choice of a constant (state-independent) friction matrix that maximizes the spectral gap of underdamped Langevin dynamics; the friction in that article does not vary with position. By contrast, the work of Lim and Tao \cite{LimTao2025} constructs a position-dependent friction tensor built from the Hessian $\nabla^2 U$, motivated by exact preconditioning of harmonic modes.   

A further line of research, distinct from those above, treats the friction as a dynamical variable governed by an auxiliary equation rather than as a chosen function of $x$.  The projective thermostatting technique of \cite{JiaLeimkuhler} applies a Nos\'e--Hoover-type thermostat through a projection onto a selected set of modes, this is particularly useful in multiscale settings where only a subset of modes require explicit temperature control. While related to the current work, the aim there was not to improve sampling, but to treat out-of-equilibrium physical models, and the method was in fact purely deterministic. Also close in spirit are the recent projection-based optimisation schemes discussed in \cite{KaLeSt2023,karoni2026adaptive}. 

There has also been a substantial body of work that replaces the underdamped SDE by a preconditioned variant in which a position-dependent matrix $B(x)$ multiplies the drift and diffusion, in close analogy with Riemannian-manifold Hamiltonian Monte Carlo \cite{GirolamiCalderhead}. A common structural difficulty of all such schemes is that preservation of the canonical measure requires the SDE drift to incorporate a divergence correction.  An analogous correction is needed for the kinetic-energy preconditioner in the underdamped formulation and computing the matrix divergence is a real impediment in practice.  Analytic preconditioners cost an additional $d$ matrix-derivative entries per step.

The ensemble preconditioner of \cite{LMW} attacks the same problem from a different angle: instead of computing $\nabla\!\cdot\!B$, the authors build $B$ from a localized empirical covariance estimated over those walkers within a distance-based cutoff of the current configuration.  The effective $B$ is then not a smooth matrix field but a discontinuous function of the joint ensemble state, depending on which walkers happen to fall inside the cutoff and on the cutoff radius itself.  Divergence identities therefore do not directly apply, and the question of measure preservation is replaced by a more delicate one about the multi-walker joint dynamics, with the cutoff radius as an additional and somewhat ad hoc hyperparameter.

The friction tensor \eqref{eq:PAL-friction} avoids the divergence problem altogether, by virtue of its incorporation into an Ornstein-Uhlenbeck part of a kinetic Langevin dynamics formulation rather than via a coordinate transformation of the physical variables.  The fluctuation--dissipation operator
$$\mathcal{L}_{\mathrm{FD}} = -(M^{-1}p)^T \Gamma(x)\nabla_p
+ \kT\,\Gamma(x):\nabla_p^2$$ preserves the canonical density
$\pi\propto e^{-H/\kT}$ pointwise for any smooth
positive-definite $\Gamma(x)$, with no drift correction involving $\nabla\!\cdot\!\Gamma$ required (Section~\ref{sec:pal}).  

The lack of any need for a divergence correction is of course understood within the Hessian-based friction of \cite{LimTao2025}.  What is specific to \eqref{eq:PAL-friction} is the specific form of $\Gamma$ itself, involving a projection matrix based on the gradient rather than the Hessian which concentrates all position dependence in the rotation of a single eigendirection $\hf(x)$.  The two eigenvalues $\alpha$ and $\alpha+\beta$ are spatially invariant, no derivative of $\Gamma$ appears anywhere in the integrator, and the canonical distribution is the correct invariant measure. For rigorous treatment of the subsequent dynamics, the hypoelliptic property and the consequent ergodicity of Langevin system with general position-dependent friction is treated in \cite{SachsLeimkuhlerDanos}; that framework covers the PAL form \eqref{eq:PAL-friction} as a special case, providing existence of an invariant measure and convergence guarantees.

\paragraph{Organisation.}
Section~\ref{sec:samadams} summarizes the SamAdams adaptive-timestep device, including the relaxation ODE governing the running sensor, options for monitor functions, choice of the Sundman transformation and the time-reweighting identity by which canonical expectations are recovered from chains generated in the rescaled time. Section~\ref{sec:pal} develops the position-adaptive Langevin (PAL) dynamics, including the specific rank-one-plus-scalar form $\Gamma(x)=\alpha I + \beta\,\hf(x)\hf(x)^T$, discuss its physical interpretation, and derives an exact $O(d)$ closed-form Ornstein--Uhlenbeck update based on it that avoids matrix decomposition. Section~\ref{app:algorithm} assembles the SamAdams scheme with the PAL friction into the AZBOBZA palindromic integrator: it defines the four elementary flows (A, Z, B, O) and explains why a particular ordering is needed. Section~\ref{sec:numerics} reports the numerical experiments on three benchmark toy model potentials (each in two spatial coordinates) and discusses parameter choices. Numerical experiments are used to compute autocorrelation functions and integrated autocorrelation times at matched gradient budget. We discuss the stepsize distribution, equilibrium accuracy of certain moments and provide an ablation study  that clarifies the contribution of the PAL friction on top of SamAdams.  Finally in Section~\ref{sec:conclusions} we summarize the findings and  outline future work, including bias analysis of the SamAdams device and extensions to higher-dimensional Bayesian and machine-learning targets.

\section{SamAdams variable timestep}
\label{sec:samadams}

In this section, we describe the variable timestepping approach of \cite{samadams} on which the SA-PAL scheme is based. Consider a system with potential $U$ that has both stiff and soft directions. If the curvature in the stiff direction is $\lambda_{\max}$, then stability of any Hamiltonian-based splitting integrator requires $h < C/\sqrt{\lambda_{\max}}$ for some method-dependent constant $C$, see for example \cite{bou2018geometric}. In a region where the curvature is mild, however, the dynamics would be stable and accurate with significantly larger steps.  A variable stepsize scheme can exploit this by automatically detecting local stiffness.

A local stiffness indicator, the monitor function, is a positive scalar function $g(x,p)$ that grows in regions where the dynamics needs to be slowed down.  Two specific choices have proved useful in the experiments of this paper:
\begin{subequations}
\label{eq:monitor}
\begin{align}
  g(x)   \;&=\; |\gU(x)|/\Omega, \qquad &&\text{(force monitor)} \label{eq:monitor-force}\\
  g(x,p) \;&=\; \bigl(|\gU(x)|^2 + |M^{-1}p|^2\bigr)/\Omega^2, \qquad &&\text{(arc-length monitor)} \label{eq:monitor-mom} 
\end{align}
\end{subequations}
The force monitor \eqref{eq:monitor-force} grows when the force is large; the scale $\Omega>0$ is chosen near the typical magnitude of $|\gU|$ over the target distribution.  This is the monitor used in the original SamAdams algorithm \cite{samadams}, where the Langevin dynamics is coupled to an auxiliary scalar relaxation that senses stiffness through $g$ and rescales the step accordingly. 

The arc-length monitor \eqref{eq:monitor-mom} is the squared norm of the Hamiltonian part of the vector field, $(\dot x,\dot p)=(M^{-1}p,-\gU)$. To ensure the stepsize does not get too large or small, we use a function $\varpsi$ and stepsizes  $\eta=\varpsi(g)$, where $\varpsi$ is decreasing in $g$. The effective step then contracts wherever this velocity is large, so the integrator advances in near-equal increments of arc length along the trajectory rather than in equal time. This results in a stochastic counterpart of the arc-length time-reparameterisation underlying adaptive integrators for deterministic Hamiltonian systems \cite{HuangLeimkuhler1997,JoLe2011}. Retaining the force term lets the step also contract in stiff regions even where the velocity is momentarily small.  This control is particularly useful near entropic bottlenecks, where the force can be small yet the geometry is tight and the velocity large.  In both cases $\Omega$ provides the global sensitivity scale and a larger $g$ signals the need for a smaller effective step.  Only the force monitor was studied in \cite{samadams}.

The SamAdams algorithm uses a relaxation approach, in particular, it introduces a sensor variable $Z$ obeying
\begin{equation}
  \frac{dZ}{dt} \;=\; -\theta\bigl(Z - g(x,p)\bigr),
  \label{eq:sensor-ode}
\end{equation}
so that $Z$ relaxes towards the current value of the monitor at rate $\theta > 0$.  A large $\theta$ makes $Z$ track $g$ closely; a small $\theta$ gives a time-averaged memory and smoother adaptation.  The averaging also helps to improve the robustness of the method in the presence of gradient noise.

For fixed $x,p$, \eqref{eq:sensor-ode} can be solved over a half-step of length $h/2$ resulting in the SamAdams update:
\begin{equation}
  Z \;\leftarrow\; e^{-\theta h/2}\,Z \;+\; (1-e^{-\theta h/2})\,g(x,p),
  \label{eq:sensor-update}
\end{equation}
which is a moving average.  The effective scale factor, $\eta$, is then defined through a Sundman kernel $\varpsi$ (see \cite{sundman1913memoire}) such that $\eta \;=\; \varpsi(Z)$, and the nominal stepsize $h$ is replaced by the effective stepsize $h_\mathrm{eff} = h\,\eta$.  We use two kernels, the first of which is defined with the force monitor; the bounded SamAdams kernel \cite{samadams},
\begin{equation}
  \varpsi^{(2)}(Z) \;=\; m_{\mathrm{lo}}\,\frac{Z^{r} + m_{\mathrm{hi}}/m_{\mathrm{lo}}}{Z^{r} + 1}
  \;\in\; [m_{\mathrm{lo}},\, m_{\mathrm{hi}}],
  \qquad 0 < m_{\mathrm{lo}} \le m_{\mathrm{hi}} \le 1,\;\; r > 0,
  \label{eq:psi-bounded}
\end{equation}
with bounds $\varpsi^{(2)}(0)=m_{\mathrm{hi}}$ and $\varpsi^{(2)}(Z{\to}\infty)=m_{\mathrm{lo}}$, and a transition centered on $Z=1$ (i.e.\ on $|\gU|=\Omega$ at equilibrium) of sharpness controlled by $r$.  The parameter $m_{\mathrm{hi}}$ caps the maximum step in quiescent regions ($m_{\mathrm{hi}} h$ should respect the stability threshold), while $m_{\mathrm{lo}}$ guarantees the chain never freezes, i.e. $h\eta \ge m_{\mathrm{lo}} h>0$ always. 

For the ``arc-length'' monitor we use instead the simple unbounded form
\begin{equation}
  \varpsi(Z) \;=\; (1+Z)^{-1/2} \;\in\;(0,1],
  \label{eq:psi-classic}
\end{equation}
which decays smoothly as $Z$ grows and imposes no upper bound, so here the nominal $h$ itself plays the role of the maximum step.  In either case $\eta\le 1$, so $h$ need only satisfy a stability condition ``on average'' and may be chosen far larger than the fixed-step threshold.

\begin{remark}
The sensor is initialized at the value of the monitor at the starting state, $Z_0 = g(x_0,p_0)$, giving $\eta_0=\varpsi(Z_0)$ immediately. The moving average introduces a warm-up transient of length $\sim 1/\theta$, where burn-in discards these steps.
\end{remark}

\begin{remark}
The time-rescaling encoded by $\varpsi$ alters the rate at which states are visited along the chain, so the canonical expectation of an observable $\phi$ is recovered by reweighting with $\eta$:
\[
  \mathbb{E}_\pi[\phi] \;\approx\;
  \frac{\sum_k \phi(x_k,p_k)\,\eta_k}{\sum_k \eta_k},
\]
which is the time-reweighted ergodic identity for the Sundman transformation (\cite{samadams}, eq.~15).  The values $\eta_k$ are stored alongside the chain and used in every reported observable; we verify all moments against reference distributions and against exact quadrature in Section~\ref{sec:numerics}.
\end{remark}

\section{Position-adaptive Langevin dynamics}
\label{sec:pal}

The underdamped Langevin dynamics with position-dependent friction and constant symmetric positive-definite mass matrix $M\in\mathrm{Sym}^{d\times d}_{+}$ are given by
\begin{align}
  dx &= M^{-1}p\,dt, \\
  dp &= -\gU(x)\,dt-\Gamma(x)M^{-1}p\,dt
        +\sqrt{2\kT\,\Gamma(x)}\,dW_t,
  \label{eq:pal-sde}
\end{align}
a minor modification of \eqref{eq:UL-x}-\eqref{eq:UL-p} for  any smooth, uniformly positive-definite friction tensor $\Gamma:\R^d\to\mathrm{Sym}^{d\times d}_{+}$. Under the usual regularity and confinement assumptions, these dynamics preserve the canonical distribution with density
\[
  \pi(x,p)
  = Z^{-1}\exp\left(-\frac{H(x,p)}{\kT}\right),
  \qquad
  H(x,p)=U(x)+\frac12 p^T M^{-1}p,
\]
for a normalisation constant $Z > 0$.

To verify invariance, consider the fluctuation--dissipation part of the forward Kolmogorov operator. Acting on a density $\rho$, it is
\[
  \mathcal{L}_{\mathrm{FD}}^*\rho
  =
  \nabla_p\cdot
  \left[
    \Gamma(x)
    \left(
      M^{-1}p\,\rho+\kT\nabla_p\rho
    \right)
  \right].
\]
Since
$
  \nabla_p\pi
  =-(\kT M)^{-1}p\,\pi,
$
it follows immediately that $\mathcal{L}_{\mathrm{FD}}^*\pi=0$ pointwise. In particular, no derivative or divergence correction involving $\Gamma(x)$ is required and the diffusion acts only in the momentum variables, while $x$ is fixed under the fluctuation--dissipation dynamics. The Hamiltonian part also preserves $\pi$, since its flow is divergence-free and conserves $H$. Hence the full forward operator satisfies $\mathcal{L}^*\pi=0$.

Under suitable growth conditions on $U$, together with uniform non-degeneracy and sufficient regularity of $\Gamma$, the resulting process is hypoelliptic and ergodic with respect to $\pi$; see \cite{SachsLeimkuhlerDanos}. Although the choice $M=I$ is often made for simplicity, a non-trivial constant mass matrix may be used to precondition anisotropic target distributions.

\subsection{The normalized rank-one friction}
We adopt and make use of the following choice of position-dependent friction
\begin{equation}
  \Gamma(x)
  =
  \alpha I+\beta\,\hf(x)\hf(x)^T,
  \qquad
  \hf(x)=\frac{\gU(x)}{|\gU(x)|},
  \label{eq:Gamma-norm}
\end{equation}
with $\alpha>0$ and $\alpha+\beta>0$\footnote{To avoid degeneracy near saddle points or extrema, i.e., when $\nabla U \approx 0$ one can instead consider 
\begin{equation*}
  \hf(x)=\frac{\gU(x)}{\sqrt{|\gU(x)|^2 + \epsilon^2}},
\end{equation*}
where $\epsilon > 0$ is small.
}. 
Since $\hf(x)\hf(x)^T$ is the orthogonal projector onto the force direction, $\Gamma(x)$ has eigenvalue $\alpha+\beta$ along $\hf(x)$ and eigenvalue $\alpha$, with multiplicity $d-1$, in the orthogonal directions. Thus, the anisotropy ratio $(\alpha+\beta)/\alpha$ is spatially constant, while the principal dissipation direction rotates with the force.

For $\beta>0$, momentum parallel to the force is damped more strongly than momentum tangent to the local level sets of $U$. This promotes rapid thermalisation in the energy-changing direction while retaining weaker damping along equipotential contours.

\subsection{Exact \texorpdfstring{$O(d)$}{O(d)} Ornstein--Uhlenbeck update}
\label{sec:ou}
For $\Gamma$ of the form \eqref{eq:Gamma-norm} with $\hat{u}$ a fixed unit vector (to denote $\hat{f}(x)$ at position $x \in \mathbb{R}^{d}$) and isotropic mass $M = m I$, the OU process given by the solution to $$dp = -\Gamma (p/m)\,dt + \sqrt{2\kT\Gamma}dW_t,$$ can be solved exactly (in the weak sense) over time $\Delta t > 0$. We then use the following convenient decomposition of the momentum $p = p_\parallel \hat{u} + p_\perp$, with $p_\parallel = \hat{u}\cdot p$, $p_\perp = p - p_\parallel\hat{u}$ whose solutions satisfy
\begin{equation}
  p_\parallel(\Delta t) = e^{-(\alpha+\beta)\Delta t/m}\,p_\parallel(0)
    + \sqrt{m\kT\bigl(1 - e^{-2(\alpha+\beta)\Delta t/m}\bigr)}\,\xi_\parallel,
  \label{eq:OU-par}
\end{equation}
\begin{equation}
  p_\perp(\Delta t) = e^{-\alpha\Delta t/m}\,p_\perp(0)
    + \sqrt{m\kT\bigl(1 - e^{-2\alpha\Delta t/m}\bigr)}\,\xi_\perp,
  \label{eq:OU-perp}
\end{equation}
where $\xi_\parallel \sim \mathcal{N}(0,1)$ and $\xi_\perp \sim \mathcal{N}(0,I_{d-1})$ are independent. Due to the fact that $\alpha/m$ and $(\alpha+\beta)/m$ are global constants, the four scalar factors in \eqref{eq:OU-par}--\eqref{eq:OU-perp} can be precomputed once for a given parameter selection $(\alpha,\beta,m,\Delta t)$; only the direction $\hat{u}$ changes at each step.

For a general mass matrix $M$, the OU update can still be performed exactly in the basis of $M^{-1}\Gamma$; see Appendix~\ref{app:mass-matrix} for a general construction.

\subsection{The {AZBOBZA} integrator}
\label{app:algorithm}
We combine the Langevin system \eqref{eq:UL-x}--\eqref{eq:UL-p} with the PAL friction \eqref{eq:Gamma-norm} and the SamAdams sensor \eqref{eq:sensor-update} into four elementary flows, each of which can be executed exactly.  The A step which advances the position using the momentum, the B step which advances the momentum using the force, the O-step as an exact Ornstein-Uhlenbeck solve as described in Section \ref{sec:ou}, all of which are present in classical ABO-based splitting integrators \cite{LeMa2013}. Then additionally, we denote the Z-step which is used for adaptive stepsizing, which given $\nu = e^{-\theta h/2}$ and the current monitor value $g(x,p)$ (either \eqref{eq:monitor-force} or \eqref{eq:monitor-mom}) is defined via \eqref{eq:sensor-update}, where $\varpsi$ is either the unbounded kernel \eqref{eq:psi-classic} or the bounded kernel \eqref{eq:psi-bounded} as appropriate.  The Z step uses the force $F=-\gU(x)$ already computed and generates no additional gradient evaluation.

We then use the splitting method ordering AZBOBZA, where the position $x$ is advanced only by the two outer $A$ half-steps; the intervening $Z$, $B$, $O$, $B$, $Z$ operators all leave $x$ fixed.  A single force $F=-\gU(x)$, evaluated once after the first $A$ half-step, therefore serves both $B$ momentum kicks, the $O$ update (through $\hat F=F/|F|$), and both $Z$ sensor updates.  A comparison of integrated autocorrelation times per step is thus a comparison per force evaluation.  As in \cite{samadams} one can show that under appropriate assumptions the AZBOBZA integrator has weak error of order two by using \cite[Proposition 6.1]{vilmart2014weak}.

\begin{remark}
The block ordering is not unique: ABZOZBA and ZOBABOZ, for instance, admit the same single-force implementation.
\end{remark}

The complete single-step update for AZBOBZA, starting from state $(x,p,Z,\eta)$ with the pre-computed constant $\nu=e^{-\theta h/2}$, is given as pseudocode as follows.
\begin{center}
\begin{tabular}{l}
\hline\\[-6pt]
\textbf{{AZBOBZA} step}\\[2pt]
\hline\\[-6pt]
1.\quad $x \leftarrow x + \tfrac{h}{2}\eta\,p$
  \hfill\textit{(A: position half-step)}\\[2pt]
2.\quad $F \leftarrow -\gU(x)$
  \hfill\textit{(single force evaluation)}\\[2pt]
3.\quad $Z \leftarrow \nu\,Z + (1{-}\nu)\,g(x,p);\quad
        \eta \leftarrow \varpsi(Z)$
  \hfill\textit{(Z: sensor half-step)}\\[2pt]
4.\quad $p \leftarrow p + \tfrac{h}{2}\eta\,F$
  \hfill\textit{(B: momentum half-step)}\\[2pt]
5.\quad $\hat{u} \leftarrow -F/|F|$;\quad
        $\Delta t \leftarrow h\eta$\\
\quad\quad $p_\parallel \leftarrow \hat{u}\cdot p$;\quad
        $p_\perp \leftarrow p - p_\parallel \hat{u}$\\
\quad\quad $p_\parallel \leftarrow e^{-(\alpha+\beta)\Delta t}\,p_\parallel
           + \sqrt{\kT(1-e^{-2(\alpha+\beta)\Delta t})}\,\xi_\parallel$\\
\quad\quad $p_\perp \leftarrow e^{-\alpha\Delta t}\,p_\perp
           + \sqrt{\kT(1-e^{-2\alpha\Delta t})}\,\xi_\perp$\\
\quad\quad $p \leftarrow p_\parallel \hat{u} + p_\perp$
  \hfill\textit{(O: PAL OU step)}\\[2pt]
6.\quad $p \leftarrow p + \tfrac{h}{2}\eta\,F$
  \hfill\textit{(B: momentum half-step, same $F$)}\\[2pt]
7.\quad $Z \leftarrow \nu\,Z + (1{-}\nu)\,g(x,p);\quad
        \eta \leftarrow \varpsi(Z)$
  \hfill\textit{(Z: sensor half-step, updated $p$)}\\[2pt]
8.\quad $x \leftarrow x + \tfrac{h}{2}\eta\,p$
  \hfill\textit{(A: position half-step)}\\[3pt]
\hline
\end{tabular}
\end{center}

\section{Numerical experiments}
\label{sec:numerics}

Before turning to a problem in higher dimensions, we test the method on three two-dimensional model potentials, chosen to exercise distinct pathologies: curved metastable basins (Rosenbrock), strongly anisotropic geometry with an entropic bottleneck (thin channel), and multiple metastable states with rare inter-well transitions at elevated temperature (Mueller--Brown). 

For each model we obtain reference marginal densities by direct numerical quadrature of the Gibbs distribution with density $\pi(x_1,x_2)\propto\exp(-U(x_1,x_2)/\kT)$ rather than by simulation.

\subsection{Model Problems}
\subsubsection{Rosenbrock (\texorpdfstring{$\kT = 1$}{kT = 1}).}
For $x_{1}$, $x_{2} \in \mathbb{R}$, we define our Rosenbrock potential function pointwise in $\mathbb{R}^{2}$ via
\[
  U(x_1, x_2) = (1-x_1)^2 + 100(x_2 - x_1^2)^2.
\]
The distribution is a narrow ridge following the Rosenbrock valley. Fixed friction and stepsize integrators such as BAOAB and OBABO typically require a moderately large friction to suppress oscillations along the ridge.

\subsubsection{Thin entropic channel (\texorpdfstring{$\kT = 1$}{kT = 1}).}
For $x_{1}$, $x_{2} \in \mathbb{R}$, we define our thin entropic channel potential function pointwise in $\mathbb{R}^{2}$ via
\[
  U(x_1, x_2) = \frac{100 x_2^2}{1 + 10 x_1^4} + 0.001(x_1^2 - 9)^2.
\]
This is constructed such that the channel narrows dramatically near $x_1 = 0$, where the transverse curvature is $200$. Fixed stepsize Langevin integrators, like BAOAB or OBABO are slow to traverse the narrow neck, often failing to pass through at all when the stepsize is large.

\subsubsection{Mueller--Brown (\texorpdfstring{$\kT = 15$}{kT = 15}).}
The Mueller--Brown potential function from \cite{MuellerBrown1979} is defined for $x_{1}$, $x_{2} \in \mathbb{R}$ by the sum of four shifted anisotropic Gaussians
\[
  U(x_1,x_2)
  \;=\; \sum_{k=1}^{4} A_k
       \exp\!\Bigl(a_k(x_1-x_{1,k})^2 + b_k(x_1-x_{1,k})(x_2-x_{2,k})
                    + c_k(x_2-x_{2,k})^2\Bigr),
\]
with the commonly used parameters
\[
\begin{array}{r|rrrr}
 k          & 1     & 2     & 3      & 4    \\ \hline
 A_k        & -200  & -100  & -170   &  15  \\
 a_k        &   -1  &   -1  &  -6.5  &  0.7 \\
 b_k        &    0  &    0  &  11    &  0.6 \\
 c_k        &  -10  &  -10  &  -6.5  &  0.7 \\
 x_{1,k}    &    1  &    0  &  -0.5  & -1.0 \\
 x_{2,k}    &    0  &  0.5  &   1.5  &  1.0
\end{array}
\]
This surface has three local minima connected by two saddle points; at $\kT = 15$ inter-well transitions occur at a physically meaningful rate.

\subsection{Parameter choices}
\label{subsec:parameters}

There are a collection of parameters both in the fixed stepsize, fixed friction integrators and the adaptive methods. We detail our choices of parameters below.
\begin{center}
\small
\begin{tabular}{lcccccc}
\toprule
Configuration & default $\gamma$ & default $h_{\mathrm{BAOAB}}$ & default $h_{\mathrm{OBABO}}$
      & $h_{\mathrm{sa}}$ & $\alpha$ & $\beta$ \\
\midrule
Rosenbrock (force)              & 5.0  & 0.02 & 0.02 & 0.8   & 0.05 & 5.0  \\
Rosenbrock (arc-length)         & 5.0  & 0.02 & 0.02 & 0.040 & 0.2  & 5.0  \\
Channel                         & 5.0  & 0.01 & 0.01 & 0.5   & 0.05 & 5.0  \\
Mueller--Brown                  & 10.0 & 0.02 & 0.02 & 0.3   & 0.5  & 10.0 \\
Mueller--Brown (arc-length)     & 10.0 & 0.02 & 0.02 & 0.1   & 0.5  & 10.0 \\
\bottomrule
\end{tabular}
\end{center}
The {AZBOBZA} nominal stepsize $h_{\mathrm{sa}}$ is chosen significantly larger than the stability-limited BAOAB/OBABO stepsizes; the SamAdams sensor automatically reduces the effective step in stiff regions. PAL parameters $(\alpha,\beta)$ are selected via a short two-stage grid sweep (described in the supplementary material); the Mueller--Brown values reflect the choice $\gamma_{\mathrm{opt}} \approx 10$ for Mueller--Brown with $\kT = 15$. The Rosenbrock force-monitor row is retained for the monitor comparison in Table~\ref{tab:rosenbrock}, whereas the arc-length row is the Rosenbrock configuration used in the main Rosenbrock figures and in the SA-only versus SA-PAL ablation.  When a figure uses a fixed-step friction or stepsize sweep, the exact fixed-step values are given in the corresponding caption or table; the entries above are the default model settings.

The SamAdams parameters (sensor relaxation rate $\theta$, kernel form, bounds $m_{\mathrm{lo}},m_{\mathrm{hi}}$, sharpness $r$, scale $\Omega$) are summarized below.  The force monitor \eqref{eq:monitor-force} with the bounded kernel \eqref{eq:psi-bounded} is used by default; the channel and the Rosenbrock arc-length configuration use the monitor \eqref{eq:monitor-mom} with the unbounded kernel \eqref{eq:psi-classic}.  On the channel this is needed because the rate-limiting feature is the entropic narrowing at $x_1=0$ (where $|\gU|\approx0$), which the force-only monitor cannot detect. The arc-length monitor is tested in the Mueller-Brown example see its effect on accuracy as well as convergence rate.
\begin{center}
\scriptsize
\begin{tabular}{lcccccc}
\toprule
Configuration & monitor & kernel & $\theta$ & $m_{\mathrm{lo}}$ & $m_{\mathrm{hi}}$ & $\Omega$ \\
\midrule
Rosenbrock (force)              & $|\gU|/\Omega$        & $\varpsi^{(2)}$ ($r=2$) & 0.1 & 0.004 & 0.025 & 20  \\
Rosenbrock (arc-length)         & $(|\gU|^2+p^TM^{-2}p)/\Omega^2$ & $(1+S)^{-1/2}$ & 0.1 & ---   & ---   & 10  \\
Channel                         & $(|\gU|^2+p^TM^{-2}p)/\Omega^2$ & $(1+S)^{-1/2}$ & 0.1 & ---   & ---   & 1.0 \\
Mueller--Brown                  & $|\gU|/\Omega$        & $\varpsi^{(2)}$ ($r=1$) & 0.2 & 0.005 & 0.200 & 100 \\
Mueller--Brown (arc-length)     & $(|\gU|^2+p^TM^{-2}p)/\Omega^2$ & $(1+S)^{-1/2}$ & 0.2 & --- & --- & 40 \\
\bottomrule
\end{tabular}
\end{center}

\subsection{Rosenbrock function}
\label{sec:rosenbrock-casestudy}
\begin{figure}[htbp]
\centering
\includegraphics[width=\linewidth]{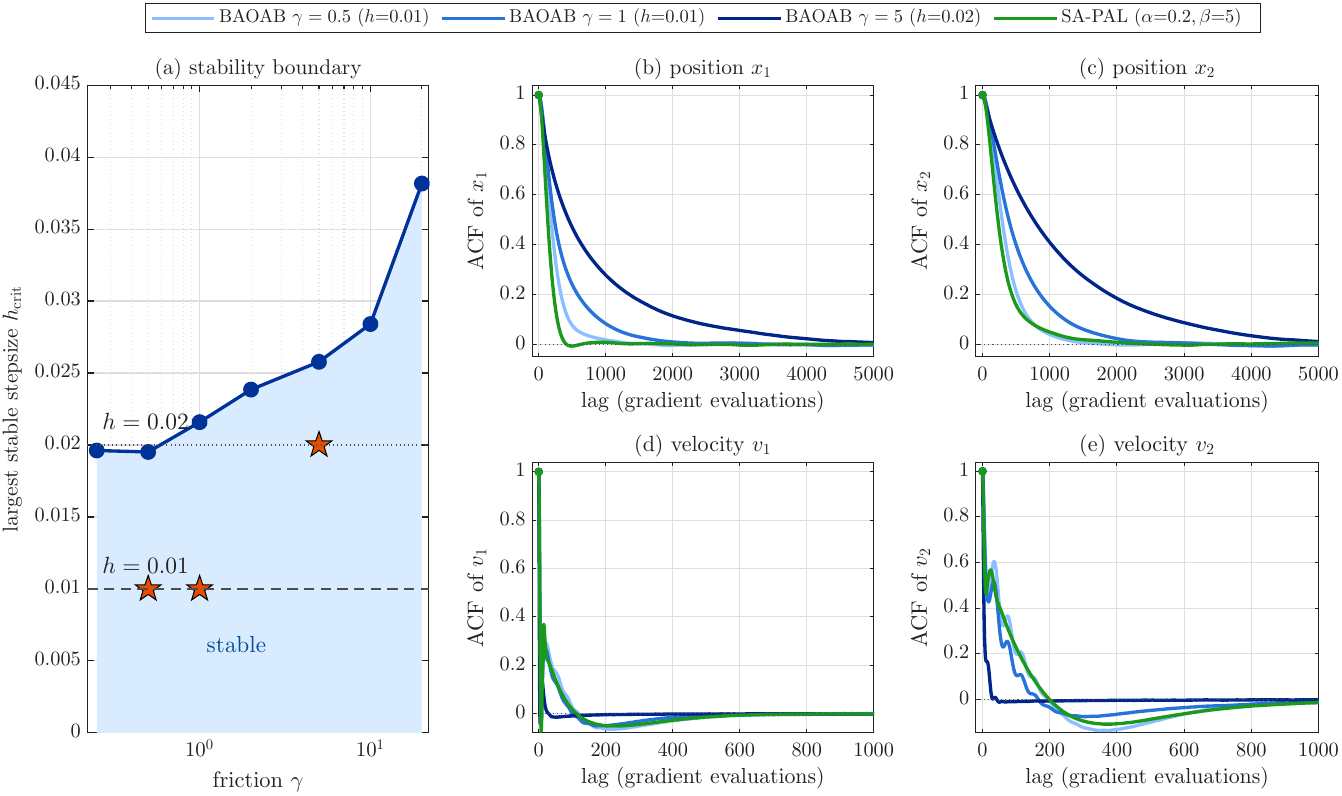}
\caption{Rosenbrock friction--stepsize dilemma. (a, left) BAOAB stability boundary $h_{\mathrm{crit}}(\gamma)$ (shaded = stable); stars mark the operating points of Table~\ref{tab:rosenbrock}. The largest stable step grows with friction, so a fast-mixing low-friction chain is forced to a small $h$. (b--e) Autocorrelation per gradient evaluation of the positions $x_1,x_2$ (top) and velocities $v_1,v_2$ (bottom) for the same four chains; curves are averaged over $9$ trajectories, three from each Rosenbrock initial condition in the code. SA-PAL (arc-length monitor, $h_{\mathrm{sa}}=0.040$, mean effective step $\approx0.0106$) decorrelates the valley coordinate faster than BAOAB at the displayed fixed-step operating points. The velocity panels are shown on a shorter lag interval because velocities decorrelate much more rapidly than positions. They show that the adaptive run remains underdamped without relying on the very heavy scalar friction that would rapidly rethermalize all momentum components.}
\label{fig:rosenbrock-dilemma}
\end{figure}

The Rosenbrock potential is a narrow, curved valley flanked by stiff walls. It exposes a dilemma intrinsic to fixed-step Langevin integrators: the largest stable stepsize depends on the friction as well as the conservative force.   At low friction the dynamics are near-ballistic, so the chain gains energy and probes the stiff walls, forcing a small stable $h$; heavier friction damps these excursions and permits a larger stepsize, but reduces the diffusion along the valley. Figure~\ref{fig:rosenbrock-dilemma}(a) shows this stability boundary $h_{\mathrm{crit}}(\gamma)$ rising with $\gamma$: $\gamma=0.5$ is unstable already at $h=0.02$ (it requires $h=0.01$), whereas $\gamma=5$ allows $h=0.02$. BAOAB is therefore squeezed, meaning that low friction values improve mixing, but limit the stepsize stability threshold. Figure~\ref{fig:rosenbrock-dilemma}(b,c) and Table~\ref{tab:rosenbrock} demonstrate the consequence: BAOAB's best achievable per-gradient IACT for the valley coordinate is about $480$, whereas the retuned arc-length SA-PAL run gives $\tau_{x_1}\approx291$.  The moderately damped stable choice $\gamma=5,h=0.02$ is nearly six times slower on $x_1$ despite taking larger fixed steps.  For the stiff transverse coordinate the arc-length SA-PAL run improves on the best low-friction fixed-step baseline and is much faster than the stable $\gamma=5,h=0.02$ row.

Rosenbrock also illustrates the difference between the two Sundman monitors of Section~\ref{sec:samadams}. With the force monitor $\varpsi^{(2)}$, which contracts the step only where $|\gU|$ is large, SA-PAL accelerates the valley coordinate $x_1$ and the energy but leaves the stiff transverse mode $x_2$ essentially at the BAOAB rate: the force is modest during the fast transverse oscillations, so a force-only monitor does not register them. The arc-length monitor \eqref{eq:monitor-mom} also responds to large velocities, and it is the configuration used in the Rosenbrock figures and ablation study.  The monitor comparison in Table~\ref{tab:rosenbrock} should therefore be read as a sensitivity check: both choices improve the slow valley coordinate substantially, while the transverse mode remains close to the best low-friction fixed-step baseline.

\begin{table}[htbp]
\centering
\caption{Rosenbrock: mixing (IACT per gradient evaluation, lower is better) and worst-case moment error for BAOAB at three frictions and for SA-PAL with each Sundman monitor. The force-monitor SA-PAL row mixes the valley coordinate $x_1$ about $3.5\times$ faster than the best BAOAB row, while the arc-length row gives a smaller but still clear improvement of about $1.7\times$. Conventions are given in Appendix~\ref{app:tables}. }
\label{tab:rosenbrock}
\begin{tabular}{lccccc}
\toprule
Method & $h$ \text{ or } $h_{\mathrm{sa}}$ & $\tau_{x_1}$ & $\tau_{x_2}$ & $\tau_{U}$ & bias (\%) \\
\midrule
BAOAB $\gamma=0.5$ & 0.01 & 482 & 726 & 351 & \textbf{0.2} \\
BAOAB $\gamma=1$   & 0.01 & 768 & 1111 & 376 & \textbf{0.2} \\
BAOAB $\gamma=5$   & 0.02 & 1673 & 2287 & 501 & 2.3 \\
SA-PAL, force monitor $\varpsi^{(2)}$       & 0.8   & \textbf{139} & 717 & \textbf{120} & 1.3 \\
{SA-PAL, arc-length monitor}         & 0.040 & 291 & 525 & 216 & 1.9 \\
\midrule
\multicolumn{6}{l}{The arc-length row uses mean effective step $N^{-1}\sum_n h_{\mathrm{sa}}\eta_n\approx0.0106$.}\\
\bottomrule
\end{tabular}
\end{table}

In light of Figure~\ref{fig:rosenbrock-dilemma}(d,e), the $\gamma=5,h=0.02$ BAOAB run is stable but visibly more damped than the low-friction BAOAB rows, and its valley-coordinate IACT remains high ($\tau_{x_1}\approx 1673$).  Running BAOAB at $\gamma=5$ with $h=0.01$ makes the physical step comparable to the adaptive run but slows the valley coordinate further.

SA-PAL keeps the isotropic friction low and uses the adaptive step to handle stiff excursions, so the velocities retain useful memory without destabilizing the trajectory.  The manuscript uses $h_{\mathrm{sa}}=0.040$, for which the ordinary mean effective step $N^{-1}\sum_n h_{\mathrm{sa}}\eta_n$ is about $0.0106$, only modestly larger than BAOAB's $h=0.01$ low-friction runs and still below the stable $\gamma=5$ fixed step. Thus the improvement is not merely a hidden larger-step comparison: it comes from combining underdamped motion with local step contraction in stiff regions.

\FloatBarrier
\subsection{The thin channel}
\label{sec:channel-casestudy}
\begin{figure}[htbp]
\centering
\includegraphics[width=\linewidth]{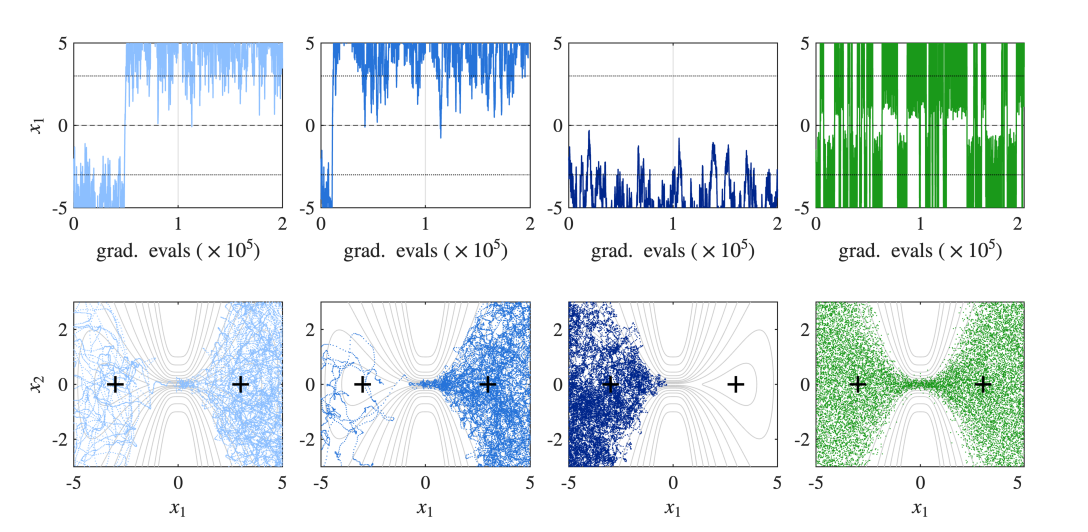}
\caption{Thin channel at equal gradient budget ($N=2\times10^5$). Columns, left to right: BAOAB at $\gamma=0.5$, $\gamma=2$, $\gamma=10$, and SA-PAL ($\alpha=0.05,\beta=5$). Top: $x_1$ vs gradient evaluations (wells at $x_1=\pm3$, neck at $x_1=0$). Bottom: trajectories on the channel potential contours. BAOAB barely crosses the neck at any friction, so its $x_1$ autocorrelation time is not reliably estimable (see Table~\ref{tab:channel}), whereas SA-PAL threads the neck repeatedly and attains $\tau_{x_1}\approx 3{,}600$.}
\label{fig:channel-traj}
\end{figure}

The thin channel isolates the regime where the SamAdams adaptive timestep is essential. The two basins at $x_1=\pm3$ communicate only through a narrow entropic neck at $x_1=0$, where the transverse confinement stiffens sharply. A fixed-step integrator must pick $h$ small enough to stay stable in the neck, yet that same $h$ is wastefully small throughout the wide basins, so each gradient evaluation advances the dynamics very little. Figure~\ref{fig:channel-traj} shows that, at an equal number of gradient evaluations, BAOAB explores poorly at every friction. In the overdamped regime ($\gamma=10$), BAOAB is trapped in one basin; in the underdamped regime, it diffuses slowly, whereas SA-PAL takes large steps in the basins and automatically shrinks $\eta$ in the neck, allowing much more frequent passes through the channel. Table~\ref{tab:channel} quantifies this: at equal cost SA-PAL makes about $16\times$ as many neck crossings as the best-tuned BAOAB and is the only scheme whose $\tau_{x_1}$ actually converges. Fixed-step BAOAB crosses the neck so rarely (a few hundred times in $2\times10^7$ steps, and even less often at the smaller $h=0.001$) that its $\tau_{x_1}$ is only a loose lower bound.

\begin{table}[htbp]
\centering
\caption{Thin channel: mixing (IACT per gradient evaluation) and the number of
neck crossings for fixed-step BAOAB versus SA-PAL at equal budget. Crossings are
the more reliable diagnostic here, since fixed-step BAOAB barely samples the slow
$x_1$ mode: SA-PAL crosses the neck far more often and is the only method to
resolve $\tau_{x_1}$. Conventions are given in Appendix~\ref{app:tables}.
$^{\dagger}$Poorly converged (too few crossings); read as a lower bound.}
\label{tab:channel}
\begin{tabular}{lccccc}
\toprule
Method & $h$ & $\tau_{x_1}$ & $\tau_{x_2}$ & $\tau_{U}$ & crossings \\
\midrule
BAOAB $\gamma=0.5$ & 0.001 & $232{,}500^\dagger$ & 67{,}900 & 41{,}600 & 33 \\
BAOAB $\gamma=0.5$ & 0.01 & $78{,}100^\dagger$ & 8{,}180 & 5{,}020 & 303 \\
BAOAB $\gamma=2$   & 0.01 & $127{,}800^\dagger$ & 28{,}100 & 16{,}600 & 278 \\
BAOAB $\gamma=10$  & 0.01 & $211{,}700^\dagger$ & 111{,}800 & 65{,}100 & 293 \\
\textbf{SA-PAL} ($\alpha{=}0.05,\beta{=}5$) & 0.5 & \textbf{3{,}602} & \textbf{247} & \textbf{206} & \textbf{4{,}857} \\
\midrule
\multicolumn{6}{l}{SA-PAL makes $\approx16\times$ more neck crossings than the best BAOAB and resolves $\tau_{x_1}$.}\\
\bottomrule
\end{tabular}
\end{table}

\subsection{Mueller--Brown}
\label{sec:mb-casestudy}

\begin{figure}[htbp]
\centering
\includegraphics[width=\linewidth]{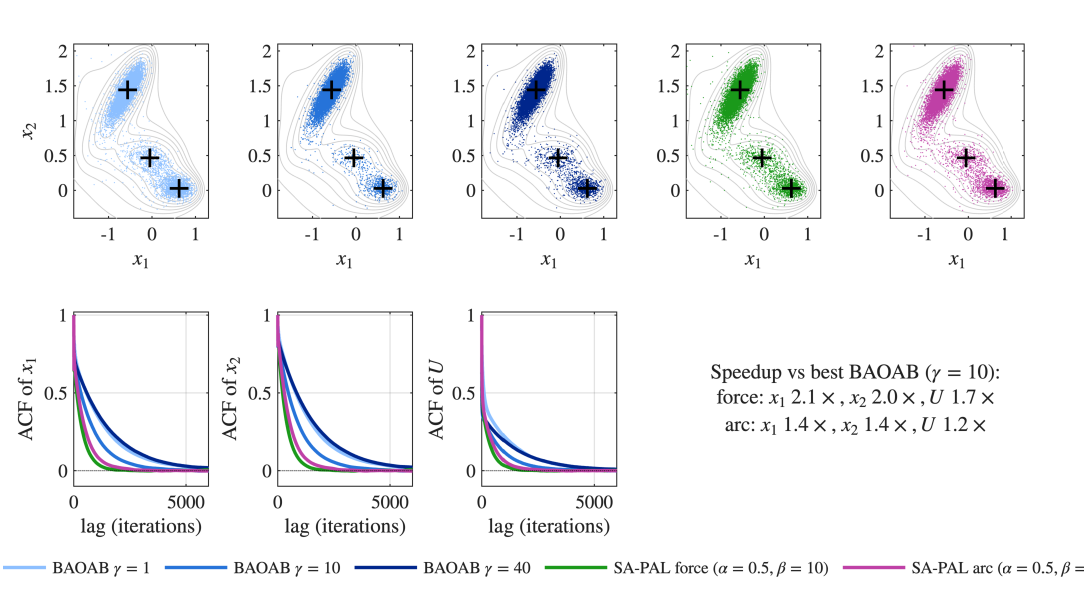}
\caption{Mueller--Brown ($kT=15$). Top: trajectories on the three-well landscape
($+$ marks the minima) at equal gradient budget, for (left to right) BAOAB at
$\gamma=1$, $\gamma=10$, $\gamma=40$, force-monitor SA-PAL, and arc-length
SA-PAL. Bottom: autocorrelation of
$x_1$, $x_2$ and $U$ (left to right) per gradient evaluation, each panel
comparing all five methods (coloured as in the top row). The force-monitor SA-PAL row decorrelates every observable fastest,
while the arc-length monitor is accurate but slower on this landscape.}
\label{fig:mb-study}
\end{figure}

The Mueller--Brown potential (with $kT=15$) is the classic three-well benchmark, with curved saddle paths connecting the basins. Here the fixed-step baseline has a genuine optimal friction: BAOAB mixes best at intermediate $\gamma\approx10$ (Table~\ref{tab:mb}), since low friction is slow to cross the saddles (under-damped, ballistic recrossing) and high friction over-damps the diffusion. Figure~\ref{fig:mb-study} (top) shows all methods explore the three wells at this temperature, while the per-variable autocorrelations (bottom) quantify the mixing: force-monitor SA-PAL ($\alpha=0.5,\beta=10$) decorrelates $x_1$, $x_2$ and $U$ faster than BAOAB at any friction. Relative to the best-tuned baseline ($\gamma=10$), the per-gradient IACT improves by about $2.1\times$ for $x_1$, $2.0\times$ for $x_2$, and $1.7\times$ for the potential energy $U$. The gain for $U$ is consistent with the structure of the dynamics: since $\dot U = \nabla U\cdot\dot x$, fluctuations of the potential energy are carried by the velocity component along the force direction $\hat f=\nabla U/\lVert\nabla U\rVert$, and SA-PAL's rank-one friction $\beta\,\hat f\hat f^{\top}$ damps precisely that component (here with $\beta=10$), driving the energy-carrying modes toward critical damping. The configurational coordinates $x_1,x_2$ instead depend on the slow inter-well transitions. The comparison is at matched accuracy: both methods exhibit a comparable discretisation bias at their operating stepsizes (about $2\%$ on the second moments, verified to be independent of chain length rather than a sampling artefact), so SA-PAL is not trading accuracy for speed. We also find the anisotropic friction is important for accuracy as well as mixing: the SamAdams-only limit $\beta=0$ samples the well populations poorly (bias $>10\%$), which the rank-one friction $\beta>0$ removes. For comparison we also include the arc-length SamAdams monitor used in the Rosenbrock and channel examples. On Mueller--Brown it is accurate but takes smaller effective steps, and consequently mixes more slowly than the force monitor.  In some sense it represents a useful compromise, offering better accuracy than  force-based SA-PAL while still giving better mixing than the accurate fixed-stepsize methods with low friction.

\begin{table}[htbp]
\centering
\caption{Mueller--Brown ($kT=15$): mixing (IACT per gradient evaluation) for
BAOAB at three frictions versus SA-PAL with force and arc-length monitors. BAOAB mixes best at intermediate friction
$\gamma=10$; the force-monitor SA-PAL row beats it on every observable at matched accuracy. The arc-length row is accurate but slower. Conventions are given in
Appendix~\ref{app:tables}.}
\label{tab:mb}
\begin{tabular}{llcccc}
\toprule
Method & monitor & $h$ & $\tau_{x_1}$ & $\tau_{x_2}$ & $\tau_{U}$ \\
\midrule
BAOAB $\gamma=1$   & fixed & 0.02 & 2092 & 2463 & 1287 \\
BAOAB $\gamma=10$  & fixed & 0.02 & 1281 & 1520 & 643 \\
BAOAB $\gamma=40$  & fixed & 0.02 & 2225 & 2617 & 1107 \\
{SA-PAL} ($\alpha{=}0.5,\beta{=}10$) & force & 0.3 & \textbf{625} & \textbf{777} & \textbf{374} \\
SA-PAL ($\alpha{=}0.5,\beta{=}10$) & arc-length & 0.1 & 902 & 1113 & 522 \\
\midrule
\multicolumn{6}{l}{Force-monitor SA-PAL vs.\ best BAOAB ($\gamma=10$): $2.1\times$ ($x_1$), $2.0\times$ ($x_2$), $1.7\times$ ($U$).}\\
\bottomrule
\end{tabular}
\end{table}

\subsection{Effective-stepsize distributions}

\begin{figure}[htbp]
\centering
\includegraphics[width=\linewidth]{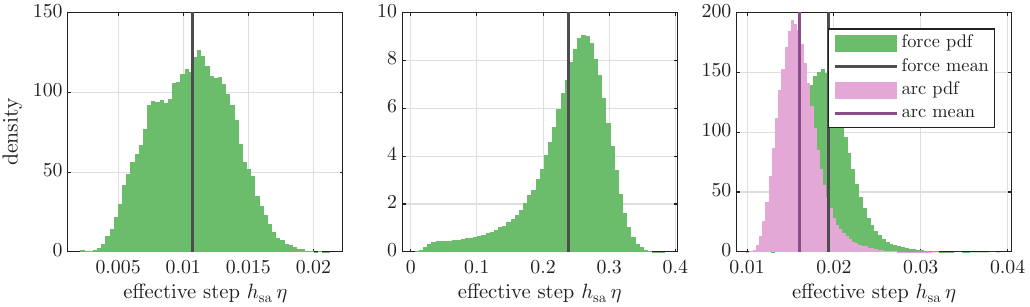}
\caption{Distribution of the SamAdams effective stepsize $h_{\mathrm{sa}}\eta$ along SA-PAL chains for Rosenbrock (left), channel (center), Mueller--Brown (right), using the tuned per-model configurations. The Mueller--Brown panel overlays the force monitor and the arc-length monitor. Solid lines mark the corresponding mean effective steps $\langle\eta\rangle h_{\mathrm{sa}}$. The spread reflects the range of stiffness encountered in each landscape ($8\times10^5$ steps).}
\label{fig:stepsize_distributions}
\end{figure}

A key diagnostic for the SamAdams device is the distribution of the effective stepsize $h_{\mathrm{sa}}\eta$ taken along converged SA-PAL chains, where $\eta = \varpsi(Z)$ is the adaptive scale factor. Figure~\ref{fig:stepsize_distributions} shows these distributions for the three tuned models, with both Mueller--Brown monitor choices overlaid in the right panel.  The arc-length kernel \eqref{eq:psi-classic} (Rosenbrock, channel, and the Mueller--Brown comparison row) produces a continuous spread, while the bounded force-monitor kernel \eqref{eq:psi-bounded} (the main Mueller--Brown row) confines $\eta$ to $[m_{\mathrm{lo}},m_{\mathrm{hi}}]$.  The contrast with the fixed BAOAB step (dashed line) is informative: on the channel the mean effective step is roughly $24\times h_{\mathrm{BAOAB}}$, exploiting the wide basins, whereas on Rosenbrock it is actually smaller than $h_{\mathrm{BAOAB}}$. SA-PAL's advantage there appears to be due to improved  in-basin mixing, not a larger step.

\subsection{The contribution of PAL on top of SamAdams}
\label{sec:ablation-pal}

The results so far combine the two devices. To isolate the contribution of PAL, we repeat the runs with $\beta=0$, which reduces the friction to the isotropic form $\gamma I$ and leaves only the SamAdams timestep active. All other settings including the nominal step, the sensor parameters, and the reweighting are unchanged. We call this variant SA-only and run it at three friction values per model. Table~\ref{tab:ablation} reports the results.

\begin{table}[h]
\centering
\caption{Ablation: SA-only ($\beta=0$, isotropic friction $\gamma I$) against SA-PAL with its tuned anisotropic friction, on all three potentials (lower is better). Removing PAL leaves the adaptive timestep alone, so the comparison isolates what PAL adds. {\color{red}Red} SA-only rows are too biased to be a fair baseline (bias $>2.5\%$); among the accuracy-matched (black) rows, bold marks the best in each column, and the ``PAL gain'' row is the fastest such SA-only IACT divided by SA-PAL's. Conventions are given in Appendix~\ref{app:tables}.}
\label{tab:ablation}
\begin{tabular}{lrrrr}
\toprule
Method & $\tau(x_1)$ & $\tau(x_2)$ & $\tau(U)$ & bias (\%) \\
\midrule
\multicolumn{5}{l}{\bfseries Rosenbrock,\quad $\kT=1$} \\
\midrule
SA-only,\; $\gamma=0.5$ & \numm{424.6}  & \numm{675.9}   & \numm{311.3} & \numm{1.98} \\
SA-only,\; {\color{red}$\gamma=2$}   & {\color{red}\numm{1239.7}} & {\color{red}\numm{1793.1}}  & {\color{red}\numm{335.0}} & {\color{red}\numm{2.51}} \\
SA-only,\; $\gamma=10$  & \numm{4576.2} & \numm{5641.2} & \numm{1113.6} & \numb{1.29} \\
{SA-PAL}\;($\alpha=0.2,\beta=5$, arc-length)
                        & \numb{288.2} & \numb{511.6} & \numb{215.8} & \numm{1.85} \\
\multicolumn{5}{l}{SA-PAL has the lowest IACTs; high isotropic friction has the smallest bias but mixes slowly.}\\
\midrule
\multicolumn{5}{l}{\bfseries Thin channel,\quad $\kT=1$} \\
\midrule
SA-only,\; {\color{red}$\gamma=0.5$} & {\color{red}\numm{4905.5}} & {\color{red}\numm{286.0}}  & {\color{red}\numm{168.6}}  & {\color{red}\numm{2.87}} \\
SA-only,\; $\gamma=2$   & \numm{7035.8} & \numm{1090.5} & \numm{611.6}  & \numb{0.56} \\
SA-only,\; $\gamma=10$  & \numm{10407.6} & \numm{3518.6} & \numm{1906.1} & \numm{1.51} \\
{SA-PAL}\;($\alpha=0.05,\beta=5$)
                        & \numb{3940.5} & \numb{245.8} & \numb{200.7} & \numm{0.69} \\
PAL gain (fastest usable SA-only $/$ SA-PAL) & 1.79 & 4.44 & 3.05 & \\
\midrule
\multicolumn{5}{l}{\bfseries M\"uller--Brown,\quad $\kT=15$} \\
\midrule
SA-only,\; {\color{red}$\gamma=1$}  & {\color{red}\numm{836.8}} & {\color{red}\numm{1096.2}} & {\color{red}\numm{565.8}} & {\color{red}\numm{6.72}} \\
SA-only,\; {\color{red}$\gamma=10$}  & {\color{red}\numm{838.0}} & {\color{red}\numm{1015.2}} & {\color{red}\numm{451.5}} & {\color{red}\numm{4.92}} \\
SA-only,\; $\gamma=40$  & \numm{1800.8} & \numm{2113.9} & \numm{886.9}  & \numb{2.11} \\
{SA-PAL force}\;($\alpha=0.5,\beta=10$)
                        & \numb{606.7} & \numb{755.6} & \numb{368.6} & \numm{2.45}\\
SA-PAL arc\;($\alpha=0.5,\beta=10$)
                        & \numm{916.1} & \numm{1129.5} & \numm{536.5} & \numm{1.22}\\
PAL gain (best usable SA-only $/$ force SA-PAL) &  2.97 &  2.80 &  2.41 &\\
\bottomrule
\end{tabular}
\end{table}

There is an important distinction between the fixed-step and adaptive-step comparisons. The fixed-step BAOAB baselines are primarily stability-limited in these examples: once the step is below the stability threshold, their equilibrium bias remains small until the step is pushed close to instability. By contrast, SamAdams and SA-PAL are often accuracy-limited before they are stability-limited. If the nominal adaptive step is made sufficiently small, both SA-only and SA-PAL produce low bias; near the largest useful adaptive steps, however, SA-PAL is noticeably more accurate than SA-only. Thus the fair comparison is not at the same aggressive nominal step when one row is biased, but at comparable sampling accuracy. Under this matched-accuracy convention, SA-PAL converges much more rapidly than the corresponding SA-only operating point.

Most of the speedup over fixed-step methods is already present in SA-only: the adaptive timestep alone accounts for the bulk of the gain. On Rosenbrock the best SA-only run is slightly faster than the best BAOAB, and on the channel it reduces the within-well autocorrelation times by an order of magnitude. PAL adds a further improvement on top when the comparison is restricted to rows with acceptable bias. On Rosenbrock this is clear-cut in mixing: SA-PAL has the lowest IACTs in Table~\ref{tab:ablation}, while the very high friction SA-only row has the smallest bias but mixes slowly. The curved valley gives the force-aligned friction a well-defined direction to act along, which is the case PAL is designed for.

Mueller--Brown is more instructive, because here an isotropic friction must cope with a trade-off: small constant friction mixes the positions quickly but samples the basins inaccurately, because the chain crosses between the wells too rarely. A large friction corrects the accuracy but slows the position mixing. Table~\ref{tab:ablation} shows that at low and intermediate friction SA-only has visibly high bias (these rows are colored red and are not used as baselines), and bringing the bias down requires a friction at which the position autocorrelation times are several times larger.

PAL with either monitor addresses this trade-off. It damps strongly along the force direction, which carries the potential energy, while leaving the slow inter-well motion lightly damped. SA-PAL therefore combines the position mixing of the low-friction case with a bias as low as the high-friction case. Against the one isotropic friction that matches its accuracy, SA-PAL mixes the slow coordinates  $2.7\times$ faster. The arc-length monitor  has better bias than the force-based SA-PAL on Mueller--Brown (as good as fixed stepsize), but its effective steps are smaller and its IACTs lie between the force-monitor SA-PAL row and the high-friction SA-only row.

\begin{figure}[htbp]
\centering
\includegraphics[width=\linewidth]{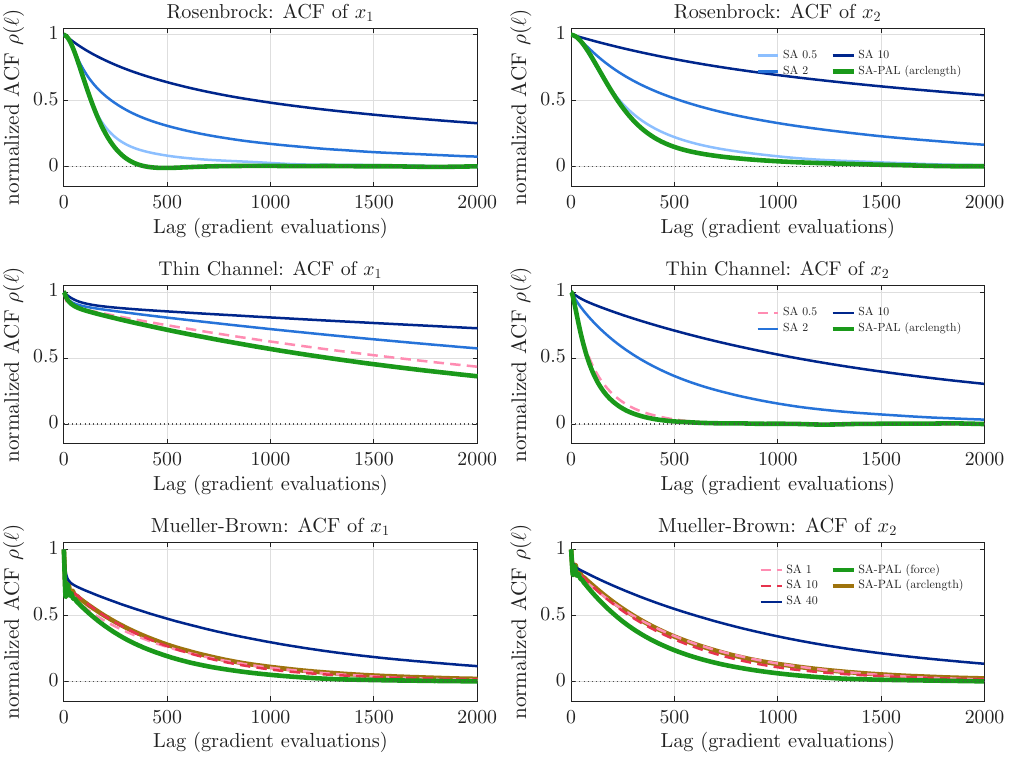}
\caption{Reweighted, normalized autocorrelation functions for the position coordinates $x_1$ and $x_2$. Each row is one model (Rosenbrock, thin channel, Mueller--Brown), and each column is one observable. The horizontal axis is lag in gradient evaluations; the vertical axis is the normalized reweighted ACF, so all curves start at one. These curves show only the correlation shape of the ratio-estimator influence sequence; the effective IACTs in Table~\ref{tab:ablation} also include the corresponding weight-variance normalisation. The curves compare SA-only at three scalar friction values with force-monitor SA-PAL (green); the Mueller--Brown row also includes the arc-length monitor comparison. Usable SA-only curves are blue; dashed pink/red SA-only curves are high-bias operating points that must be rejected. All chains use the same nominal step $h_{\mathrm{sa}}$ and the same SamAdams sensor settings; the only change is $\beta=0$ for SA-only versus the tuned $\beta>0$ for SA-PAL. Dashed SA-only curves are friction values whose sampling bias is well above SA-PAL's (Table~\ref{tab:ablation}), so their faster decay does not correspond to a usable operating point.}
\label{fig:sa-vs-sapal}
\end{figure}

An intermediate friction does not help: it is faster but still too biased to use. For this reason the low-friction SA-only curves in Figure~\ref{fig:sa-vs-sapal} are drawn dashed in pink/red, since their faster decay does not correspond to correct sampling; the only accurate (solid) SA-only option has much slower mixing than SA-PAL. We set $\beta$ to the largest value that keeps the bias below the fixed-step reference; a smaller $\beta$ gives slightly faster position mixing at the cost of more bias.

\subsection{Equilibrium accuracy at fixed budget}
Having established the mixing efficiency on each model and the separate contribution of each device, we finally report the accuracy of ensemble-averaged observables.  All methods are run for the same number of gradient evaluations, so the comparison is at matched computational budget.

\begin{figure}[htbp]
\centering
\includegraphics[width=\linewidth]{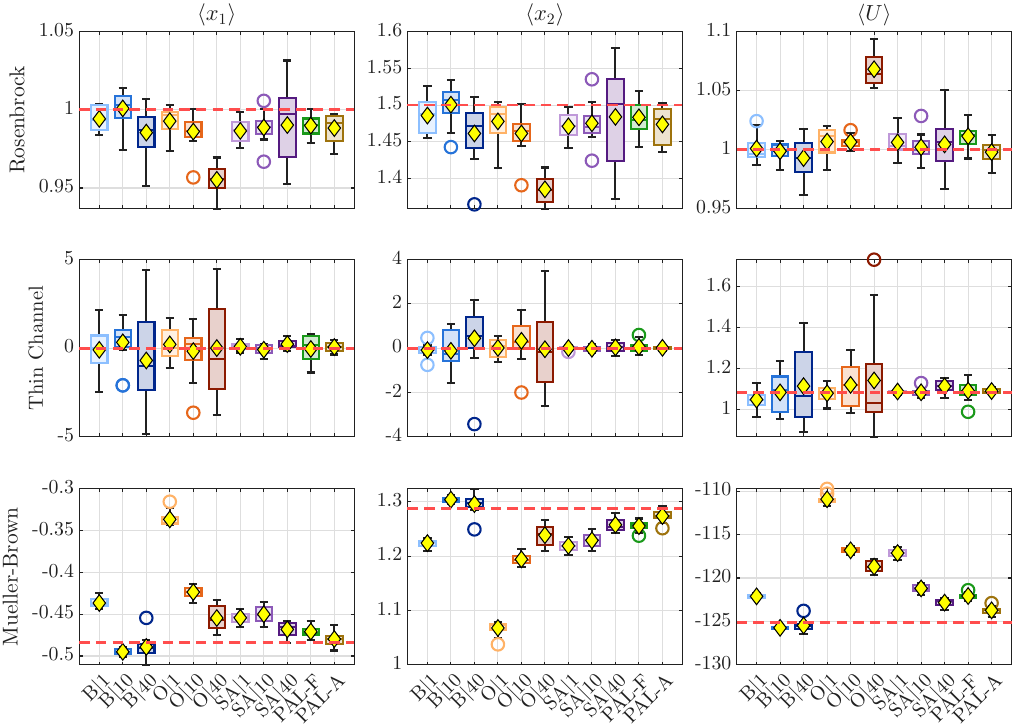}
\caption{Per-chain estimates of $\langle x_1\rangle$, $\langle x_2\rangle$ and $\langle U\rangle$ (columns, left to right) from independent BAOAB, OBABO, SA-only and SA-PAL (F=force, A= arclength) chains at matched gradient budget, plotted against the exact target values from quadrature (red dashed lines). Rows: Rosenbrock, thin channel, Mueller--Brown. BAOAB, OBABO and SA-only are each shown at three frictions; box width is the run-to-run sampling spread, box position relative to the target line is the systematic bias.}
\label{fig:accuracies}
\end{figure}

Figure~\ref{fig:accuracies} summarizes these estimates. On Rosenbrock SA-PAL is the most accurate scheme and also the least variable: its boxes lie closest to the target lines and are the narrowest of any method. In the thin channel there is little gain in the invariant measure accuracy for SA-PAL over SA-only. In the Mueller-Brown model the accuracy of SA-PAL roughly matches that of BAOAB at sufficiently high friction (of course achieving this with much less computation due to the lower IACT).   The fixed-step integrators behave as one expects when sampling the invariant measure. BAOAB is generally accurate given enough computation, and is typically most accurate at the higher frictions, that is, where it mixes slowest. Therefore its accuracy and its mixing speed pull in opposite directions. OBABO is the least accurate overall: its bias is largest at low friction, and conspicuously so on Mueller--Brown, where the low-friction chains miss the target by a wide margin. SA-only generally improves on both fixed-step methods, most clearly on the thin channel, where the adaptive step samples the slow coordinate much better at the same budget, but does not help much in the Mueller--Brown. The added Mueller--Brown arc-length box is close to the target lines, confirming that it has improved invariant-measure bias compared to SA-PAL but at the cost of some mixing efficiency.

\subsection{Bayesian linear regression with horseshoe prior}
\label{sec:horseshoe}

Our last example is a Bayesian linear regression with a horseshoe prior \cite{Carvalho2010}. This is a higher-dimensional target, and a harder one. The horseshoe is a shrinkage prior, meaning it pulls small coefficients towards zero while leaving the large ones essentially alone, and the price is a posterior that is sharply curved near the origin. We use $N=100$ observations and $K=25$ predictors, with true coefficients $\beta_{1:4} = [3, -2, 1.5, 1]$ and $\beta_{5:25} = 0$, noise scale $\sigma=0.5$, and a design matrix $X$ with i.i.d.\ Gaussian entries $X_{ij} \sim \mathcal{N}(0, 1/N)$. Writing the state as $w = [\beta, \eta, \xi]^T \in \mathbb{R}^{2K+1}$, the potential is
\begin{equation}
    U(\beta, \eta, \xi) = \frac{\|y - X\beta\|^2}{2\sigma^2} + \sum_{j=1}^K \frac{\beta_j^2}{2 e^{2(\xi + \eta_j)}} + (K-1)\xi + \sum_{j=1}^K \log(1 + e^{2\eta_j}) + \log(1 + e^{2\xi}),
\end{equation}
where $X \in \mathbb{R}^{N \times K}$ is the design matrix, $y \in \mathbb{R}^N$ is the response, and $\sigma$ is the noise scale. The difficulty lives in the $(\beta_j, \xi+\eta_j)$ plane, where the posterior has the familiar ``funnel'' shape. When a coefficient sits near zero its scale $\xi+\eta_j$ can shrink as well, leaving a narrow neck of high curvature; further out the funnel opens up and large steps are safe. A single fixed step cannot serve both ends: a step small enough for the neck is wasteful in the mouth, and a step sized for the mouth is unstable in the neck. We give SamAdams a step-size floor of $m=10^{-4}$ so the step can drop far enough to stay stable in the neck.

We compared  BAOAB, SA-only ($\beta=0$), and SA-PAL, using 20 independent seeds for each configuration, looking at three distinct fidelity regimes. For SA-PAL we used $h_{\mathrm{sa}}=0.1$ in the high-fidelity regime, $h_{\mathrm{sa}}=0.2$ in the moderate regime, and $h_{\mathrm{sa}}=0.28$ in the high-performance regime, reduced slightly from $0.3$ to keep the bias below $2.5\%$ tolerance. For SA-only we chose $h_{\mathrm{sa}}$ by the same bias rule: within each regime we swept the nominal SA-only step downward and kept the largest value whose estimated bias in $\langle|\beta|^2\rangle$ was below $2.5\%$. This gives $h_{\mathrm{sa}}=0.1$, $0.16$, and $0.20$ for the high-fidelity, moderate, and high-performance SA-only rows, respectively. BAOAB was run at $h \in \{10^{-3}, 5\cdot 10^{-3}, 10^{-2}\}$; on this funnel its largest reliably stable step is $5\cdot 10^{-3}$, while $h=10^{-2}$ already sits at the edge of stability and is seed-sensitive (some runs diverge). That fixed-step BAOAB cannot be pushed to such a step is precisely the regime the adaptive scheme is built for.

\begin{figure}[htbp]
\centering
\includegraphics[width=\linewidth]{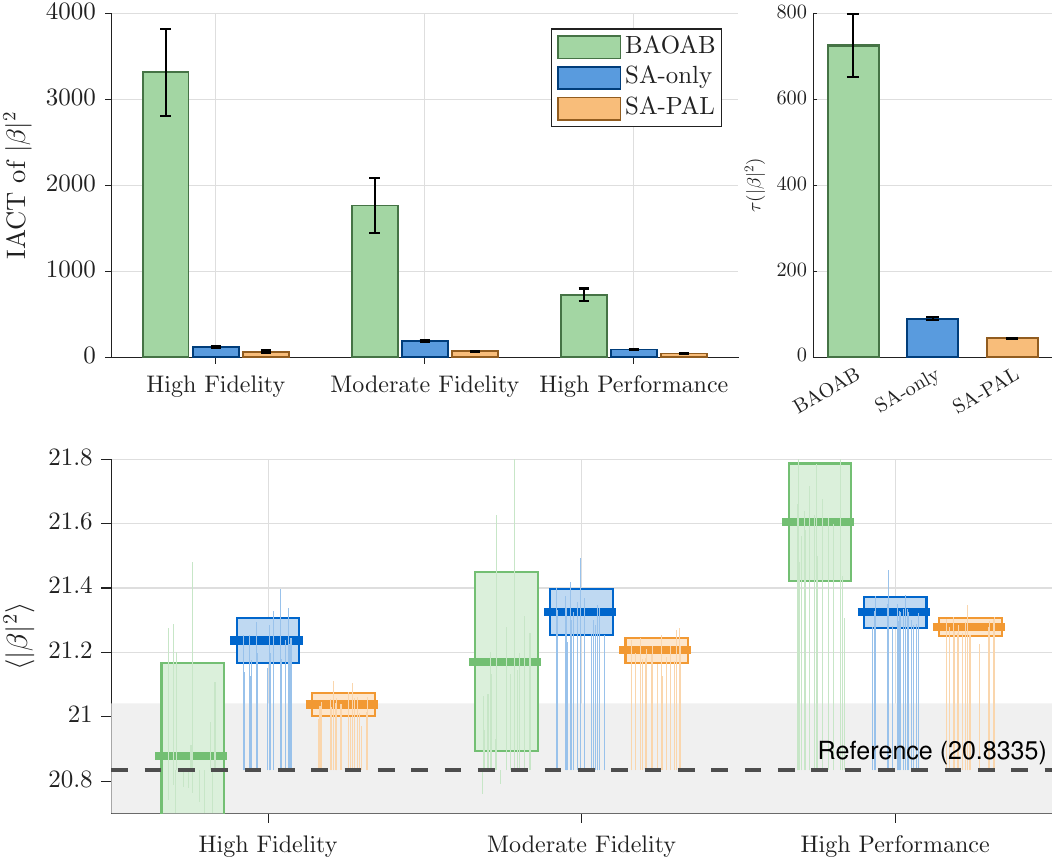}
\caption{Comparison of BAOAB, SA-only, and SA-PAL on the horseshoe regression problem across 20 independent seeds. Top: IACT of the aggregate observable $|\beta|^2$ as a grouped bar chart (bars = mean over seeds, black caps = $\pm$ standard deviation), with a companion panel on the right zooming into the high-performance regime. Bottom: sampling accuracy for the second moment $\langle|\beta|^2\rangle$, shown as box plots (box = mean $\pm$ standard deviation, thick line = mean) with thin vertical needle lines for each individual run. Among rows satisfying the $2.5\%$ bias tolerance, SA-PAL gives the best aggregate mixing.}
\label{fig:horseshoe}
\end{figure}

Figure~\ref{fig:horseshoe} collects the results. The top panel reports the effective IACT of the aggregate coefficient observable $|\beta|^2$, the same quantity used for the accuracy check in the bottom panel. For reweighted SamAdams chains this effective IACT is computed from the long-run variance of the delta-method ratio-estimator influence sequence and normalized by the weighted target variance of $|\beta|^2$, as described in Appendix~\ref{app:tables}. With this convention, SA-PAL is the fastest aggregate mixer in all three regimes among rows satisfying the $2.5\%$ bias tolerance. The gains over BAOAB are large because the fixed-step chain must remain stable in the narrow funnel neck. After bias matching the SA-only rows, SA-PAL reduces $\tau(|\beta|^2)$ by factors of about $1.88$, $2.70$, and $2.04$ in the high-fidelity, moderate, and high-performance regimes, respectively. The bottom panel shows why the bias convention matters: increasing the SA-only step beyond the selected values gives faster but biased estimates of $\langle|\beta|^2\rangle$, whereas the reported SA-PAL and SA-only rows remain close to the reference value $\langle|\beta|^2\rangle \approx 20.83$.

The numerical IACT and accuracy values corresponding to Figure~\ref{fig:horseshoe} are reported in Table~\ref{tab:horseshoe}. Here the mixing observable is $|\beta|^2$, matching the posterior moment whose accuracy is reported in the same table. Red rows have bias above $2.5\%$ and are not used as accuracy-matched baselines. With this convention, SA-PAL is the best usable aggregate mixer in each fidelity regime. The SA-only steps in the last two regimes are smaller than the nominal aggressive values precisely because they were retuned to satisfy the same bias cutoff; without this retuning SA-only appears faster but no longer gives a fair accuracy-matched comparison. Retuning the comparison for individual signal coordinates $\beta_1,\beta_2$ largely removes the SA-PAL advantage over SA-only; the claim here is therefore deliberately about the aggregate second-moment observable, and it remains sensitive to the tuning criterion.

\begin{table}[htbp]
\centering
\small
\caption{Horseshoe regression: aggregate mixing and accuracy across three
fidelity regimes. The IACT is computed for $|\beta|^2$, and the bias is the
relative error in $\langle|\beta|^2\rangle$ against the reference value
$20.8335$. The final column reports the speedup over the BAOAB row in the same
regime. SA-only uses the largest swept $h_{\mathrm{sa}}$ whose bias remains
below $2.5\%$ in that regime. {\color{red}Red} rows have bias $>2.5\%$ and
are not used as accuracy-matched operating points.}
\label{tab:horseshoe}
\begin{tabular}{@{}lrrrr@{}}
\toprule
Method & $\tau(|\beta|^2)$ & $\langle|\beta|^2\rangle$ & Bias (\%) & Speedup \\
\midrule
\multicolumn{5}{l}{\bfseries High Fidelity} \\
\midrule
BAOAB ($h=10^{-3}$)      & \numm{3315.6} & \numm{20.878} & \numb{0.22} & \numm{1.00} \\
SA-only ($h_{\mathrm{sa}}=0.10$) & \numm{119.3}  & \numm{21.237} & \numm{1.94} & \numm{27.79} \\
\textbf{SA-PAL}          & \numb{63.3} & \numm{21.038} & \numm{0.98} & \numb{52.38} \\
\midrule
\multicolumn{5}{l}{\bfseries Moderate Fidelity} \\
\midrule
BAOAB ($h=5\cdot 10^{-3}$) & \numm{1764.7} & \numm{21.171} & \numb{1.62} & \numm{1.00} \\
SA-only ($h_{\mathrm{sa}}=0.16$) & \numm{186.8} & \numm{21.325} & \numm{2.36} & \numm{9.45} \\
\textbf{SA-PAL}          & \numb{69.3}  & \numm{21.206} & \numm{1.79} & \numb{25.46} \\
\midrule
\multicolumn{5}{l}{\bfseries High Performance} \\
\midrule
{\color{red}BAOAB ($h=10^{-2}$)$^{\dagger}$} & {\color{red}\numm{725.2}} & {\color{red}\numm{21.605}} & {\color{red}\numm{3.70}} & {\color{red}\numm{1.00}} \\
SA-only ($h_{\mathrm{sa}}=0.20$) & \numm{89.2} & \numm{21.325} & \numm{2.36} & \numm{8.13} \\
\textbf{SA-PAL} ($h_{\mathrm{sa}}=0.28$) & \numb{43.7} & \numm{21.279} & \numb{2.14} & \numb{16.59} \\
\midrule
\multicolumn{5}{l}{SA-PAL/SA-only IACT ratios after bias matching: $1.88\times$, $2.70\times$, and $2.04\times$.}\\
\midrule
\multicolumn{5}{l}{$^{\dagger}$$h=10^{-2}$ is at the edge of fixed-step stability on the funnel and seed-sensitive (some}\\
\multicolumn{5}{l}{\phantom{$^{\dagger}$}runs diverge); not a robust operating point. The largest reliably stable BAOAB step is $5\cdot10^{-3}$.}\\
\bottomrule
\end{tabular}
\end{table}

\section{Conclusions}
\label{sec:conclusions}

We have introduced the {AZBOBZA} integrator, combining SamAdams adaptive timesteps with position-adaptive Langevin (PAL) dynamics. The scheme retains the $O(d)$ cost per step of  BAOAB/OBABO integrators while adapting simultaneously to local stiffness (via $\eta$) and local force direction (via the PAL friction tensor).  Numerical experiments on various distinct benchmark potentials demonstrate substantially reduced integrated autocorrelation times and accurate reproduction of the canonical distribution.

The free parameters $(\alpha, \beta, \theta, h_{\mathrm{sa}})$ have interpretable roles: $\alpha$ and $\alpha+\beta$ set the transverse and parallel friction; $\theta^{-1}$ sets the memory of the stiffness sensor; and $h_{\mathrm{sa}}$ is a nominal stepsize that is automatically reduced in stiff regions.  Nevertheless, the numerical experiments also show that the tuning problem is not entirely trivial.  Fixed-step BAOAB is usually limited by stability, whereas SamAdams-based schemes can become limited by equilibrium accuracy at large nominal steps even while the trajectories remain stable. Consequently the meaningful comparison between SA-only and SA-PAL is at matched bias: all convergent methods have low bias when the step is small, but at larger useful adaptive steps SA-PAL is more accurate than SA-only, and when SA-only is retuned to match that accuracy SA-PAL converges substantially faster.  These parameters are coupled through the discrete stability threshold, the SamAdams sensor response, and the observable used to measure mixing.  In the horseshoe example, for instance, SA-PAL gives a clear advantage for the aggregate second-moment observable, while retuning around individual signal coordinates can largely remove the advantage over SA-only.  Developing a principled, problem-independent tuning procedure for the combined SA-PAL scheme therefore remains an important open issue.

Future directions include: a rigorous analysis of the bias introduced by the SamAdams sensor; a systematic theory and practical protocol for tuning the coupled parameters; extension to higher-dimensional targets from Bayesian machine learning; and the combination of AZBOBZA with splitting methods for stochastic gradient estimators.

\appendix

\section{Langevin dynamics with a general mass matrix}
\label{app:mass-matrix}

Consider the underdamped Langevin system with a constant, symmetric positive-definite mass matrix $M$ and Hamiltonian $H(x,p) = U(x) + \tfrac{1}{2} p^T M^{-1} p$. The Ornstein--Uhlenbeck (OU) process for the momentum $p$ is given by the SDE
\begin{equation}
    dp = -\Gamma M^{-1} p \, dt + \sqrt{2\kT\Gamma} \, dW_t,
\end{equation}
which preserves the marginal canonical density $\pi_p(p) \propto \exp(-\tfrac{1}{2} p^T M^{-1} p / \kT)$ for any smooth positive-definite friction tensor $\Gamma$. To maintain cheap computational cost when $M \neq I$, we define the friction tensor in a way that is compatible with the kinetic energy structure:
\begin{equation}
    \Gamma(x) = \alpha M + \beta \hf(x) \hf(x)^T, \qquad \hf(x) = \frac{\gU(x)}{|\gU(x)|}.
    \label{eq:Gamma-mass}
\end{equation}
The exact integration of the OU process can be performed efficiently by transforming to the decorrelated variable $q = M^{-1/2} p$, in which the kinetic energy is simply $\tfrac{1}{2} |q|^2$. Precomputing $M^{1/2}$ and $M^{-1/2}$ (which is done only once), the SDE for $q$ becomes
\begin{equation}
    dq = -\mathcal{A}(x) q \, dt + \sqrt{2\kT \mathcal{A}(x)} \, dW_t, \qquad \mathcal{A}(x) = M^{-1/2} \Gamma(x) M^{-1/2}.
\end{equation}
Substituting \eqref{eq:Gamma-mass} into the expression for $\mathcal{A}$ yields:
\begin{equation}
    \mathcal{A}(x) = \alpha I + \beta \hat{v}(x) \hat{v}(x)^T, \qquad \hat{v}(x) = M^{-1/2} \hf(x).
\end{equation}
This is a rank-one update to the identity matrix. Let $w(x) = \hat{v}(x)/|\hat{v}(x)|$ be the normalized direction and $\lambda(x) = \beta |\hat{v}(x)|^2$. Then $\mathcal{A} = \alpha I + \lambda w w^T$, which has eigenvalues $\alpha$ (with multiplicity $d-1$) and $\alpha + \lambda$ (with multiplicity 1). The process for $q$ is thus an isotropic OU process with a single directional enhancement, which can be solved exactly using \eqref{eq:OU-par}--\eqref{eq:OU-perp} with parameters $(\alpha, \lambda)$ and projection vector $w$. The momentum is then recovered as $p = M^{1/2} q$.

\section{Conventions for the results tables}
\label{app:tables}
Unless stated otherwise, Tables~\ref{tab:rosenbrock}, \ref{tab:channel}, \ref{tab:mb}, \ref{tab:ablation}, and \ref{tab:horseshoe} follow the conventions below.

\paragraph{Mixing.} We report the integrated autocorrelation time (IACT) $\tau$ per gradient evaluation, estimated via an automatically truncated self-consistent automatic windowing rule \cite{SokalMCLecture} and averaged over independent seeds ($8$ seeds, or $24$ for Rosenbrock). For fixed-step chains, the ACF is the usual normalized ACF of the centered observable. For SamAdams-reweighted chains, the target estimator is the self-normalized ratio
\[
  \hat\mu_\phi=\frac{\sum_n w_n\phi_n}{\sum_n w_n}, \qquad w_n=\eta_n,
\]
so the reported ACF is computed from the ratio-estimator influence
sequence
\[
  Y_n = w_n(\phi_n-\hat\mu_\phi).
\]
This is the delta-method linearisation of the ratio estimator, due to the fact that the random denominator contributes through the centering term. The normalized ACF curves therefore show the correlation structure of $Y_n$. For the reported scalar IACT, however, we first estimate the long-run variance
\[
  \operatorname{LRV}(Y)=\operatorname{Var}(Y_n)+2\sum_{\ell\ge1} \operatorname{Cov}(Y_n,Y_{n+\ell})
\]
using the same automatically truncated window, and then report the estimator effective value
\[
  \tau_{\mathrm{eff}}(\phi) = \frac{\operatorname{LRV}(Y)}{\bar w^2\,\widehat{\operatorname{Var}}_\pi(\phi)}, \qquad \widehat{\operatorname{Var}}_\pi(\phi) = \frac{\sum_n w_n(\phi_n-\hat\mu_\phi)^2}{\sum_n w_n}.
\]
This reduces to the ordinary IACT when $w_n\equiv1$, while for reweighted chains it retains the denominator fluctuations through $Y_n$ and normalizes by the target variance of the original observable rather than by $\operatorname{Var}(Y_n)$. This convention is standard in Markov-chain error analysis for reweighted ratio estimators \cite{RoyTanFlegal2018}; the long-run variance is estimated from the autocovariance sum as in standard spectral/batch-means output analysis \cite{FlegalJones2010}. Every integrator uses one gradient per step, so IACT per gradient is a cost-normalized measure of mixing; lower is better, and bold marks the best value in each column.

\paragraph{Accuracy.} The bias column is the largest relative error, over the moments $\langle x_1\rangle, \langle x_2\rangle, \langle x_1^2\rangle$ and $\langle U\rangle$, between the sampled estimate and the exact value obtained by quadrature; first-moment errors are measured in units of the coordinate standard deviation. On the Rosenbrock and channel marginals the quartic moment $\langle x_2^2\rangle=\langle x_1^4\rangle$ is excluded, because at these chain lengths its estimate is dominated by sampling variance rather than discretisation error. In the ablation and horseshoe tables, red rows have bias above $2.5\%$ and are excluded from accuracy-matched baseline comparisons. We use this same $2.5\%$ cutoff throughout the paper when deciding whether a row is an acceptable operating point for matched-accuracy comparisons. The precise value is necessarily somewhat arbitrary, as any finite-sample bias cutoff would be. We chose it as a pragmatic compromise: small enough to rule out visibly biased aggressive parameter choices, but large enough that the required ensembles finish in a reasonable amount of computation time.

\paragraph{Stepsizes.} BAOAB and OBABO run at their stability-limited fixed step; SA-PAL runs at a larger nominal timestep whose effective value $h_{\mathrm{sa}}\langle\eta\rangle$ is automatically reduced in stiff regions. All per-model parameters are listed in Section~\ref{subsec:parameters}.

\section*{Acknowledgements}
The work described in this article was initiated while both authors were guests of the Erwin Schr\"odinger Institute (ESI) at the University of Vienna in May of 2026; the work was completed during the workshop ``Scalable MCMC Sampling'' held at the Forschungsinstitut f\"ur Mathematik (FIM) at ETH Z\"urich in June of 2026. We are greatly indebted to these two institutes for making this work possible. We also thank Christoph Dellago, Jianfeng Lu, and Gilles Vilmart for helpful discussions.


\begin{thebibliography}{24}
\providecommand{\natexlab}[1]{#1}
\providecommand{\url}[1]{\texttt{#1}}
\expandafter\ifx\csname urlstyle\endcsname\relax
  \providecommand{\doi}[1]{doi: #1}\else
  \providecommand{\doi}{doi: \begingroup \urlstyle{rm}\Url}\fi

\bibitem[Bou-Rabee and Sanz-Serna(2018)]{bou2018geometric}
N.~Bou-Rabee and J.~M. Sanz-Serna.
\newblock Geometric integrators and the {H}amiltonian {M}onte {C}arlo method.
\newblock \emph{Acta Numerica}, 27:\penalty0 113--206, 2018.
\newblock \doi{10.1017/S0962492917000101}.

\bibitem[Cao et~al.(2023)Cao, Lu, and Wang]{cao2023explicit}
Y.~Cao, J.~Lu, and L.~Wang.
\newblock On explicit {$L^2$}-convergence rate estimate for underdamped
  {L}angevin dynamics.
\newblock \emph{Archive for Rational Mechanics and Analysis}, 247\penalty0
  (5):\penalty0 90, 2023.
\newblock \doi{10.1007/s00205-023-01922-4}.

\bibitem[Carvalho et~al.(2010)Carvalho, Polson, and Scott]{Carvalho2010}
C.~M. Carvalho, N.~G. Polson, and J.~G. Scott.
\newblock The horseshoe estimator for sparse signals.
\newblock \emph{Biometrika}, 97\penalty0 (2):\penalty0 465--480, 2010.
\newblock \doi{10.1093/biomet/asq017}.

\bibitem[Chak et~al.(2023)Chak, Kantas, Leli{\`e}vre, and Pavliotis]{Chak2023}
M.~Chak, N.~Kantas, T.~Leli{\`e}vre, and G.~A. Pavliotis.
\newblock Optimal friction matrix for underdamped {L}angevin sampling.
\newblock \emph{ESAIM Math. Model. Numer. Anal.}, 57\penalty0 (6):\penalty0
  3335--3371, 2023.
\newblock \doi{10.1051/m2an/2023083}.

\bibitem[Flegal and Jones(2010)]{FlegalJones2010}
J.~M. Flegal and G.~L. Jones.
\newblock Batch means and spectral variance estimators in {M}arkov chain
  {M}onte {C}arlo.
\newblock \emph{Annals of Statistics}, 38\penalty0 (2):\penalty0 1034--1070,
  2010.
\newblock \doi{10.1214/09-AOS735}.

\bibitem[Girolami and Calderhead(2011)]{GirolamiCalderhead}
M.~Girolami and B.~Calderhead.
\newblock Riemann manifold {L}angevin and {H}amiltonian {M}onte {C}arlo
  methods.
\newblock \emph{J. R. Stat. Soc. B}, 73\penalty0 (2):\penalty0 123--214, 2011.
\newblock \doi{10.1111/j.1467-9868.2010.00765.x}.

\bibitem[Huang and Leimkuhler(1997)]{HuangLeimkuhler1997}
W.~Huang and B.~Leimkuhler.
\newblock The adaptive {V}erlet method.
\newblock \emph{SIAM Journal on Scientific Computing}, 18\penalty0
  (1):\penalty0 239--256, 1997.
\newblock \doi{10.1137/S1064827595284658}.

\bibitem[Jia and Leimkuhler(2005)]{JiaLeimkuhler}
Z.~Jia and B.~Leimkuhler.
\newblock A projective thermostatting dynamics technique.
\newblock \emph{Multiscale Model. Simul.}, 4\penalty0 (2):\penalty0 563--583,
  2005.
\newblock \doi{10.1137/040603863}.

\bibitem[Jones and Leimkuhler(2011)]{JoLe2011}
A.~Jones and B.~Leimkuhler.
\newblock Adaptive stochastic methods for sampling driven molecular systems.
\newblock \emph{Journal of Chemical Physics}, 135\penalty0 (8):\penalty0
  084125, 2011.
\newblock \doi{10.1063/1.3626941}.

\bibitem[Karoni et~al.(2023)Karoni, Leimkuhler, and Stoltz]{KaLeSt2023}
A.~Karoni, B.~Leimkuhler, and G.~Stoltz.
\newblock Friction-adaptive descent: {A} family of dynamics-based optimization
  methods.
\newblock \emph{Journal of Computational Dynamics}, 10\penalty0 (4):\penalty0
  450--484, 2023.
\newblock \doi{10.3934/jcd.2023007}.

\bibitem[Karoni et~al.(2026)Karoni, Rajpal, Leimkuhler, and
  Stoltz]{karoni2026adaptive}
A.~Karoni, R.~Rajpal, B.~Leimkuhler, and G.~Stoltz.
\newblock Adaptive momentum and nonlinear damping for neural network training.
\newblock \emph{arXiv preprint arXiv:2602.00334}, 2026.

\bibitem[Leimkuhler and Matthews(2013)]{LeMa2013}
B.~Leimkuhler and C.~Matthews.
\newblock Rational construction of stochastic numerical methods for molecular
  sampling.
\newblock \emph{Appl. Math. Res. eXpress}, 2013\penalty0 (1):\penalty0 34--56,
  2013.
\newblock \doi{10.1093/amrx/abs010}.

\bibitem[Leimkuhler et~al.(2024)Leimkuhler, Paulin, and Whalley]{LePaWh2024}
B.~Leimkuhler, D.~Paulin, and P.~A. Whalley.
\newblock Contraction and convergence rates for discretized kinetic {L}angevin
  dynamics.
\newblock \emph{SIAM Journal on Numerical Analysis}, 62\penalty0 (3):\penalty0
  1226--1258, 2024.
\newblock \doi{10.1137/23M1556289}.

\bibitem[Leimkuhler et~al.(2026)Leimkuhler, Lohmann, and Whalley]{samadams}
B.~Leimkuhler, R.~Lohmann, and P.~A. Whalley.
\newblock A {L}angevin sampling algorithm inspired by the {A}dam optimizer.
\newblock \emph{ACM Trans. Probab. Mach. Learn.}, 1, 2026.
\newblock \doi{10.1145/3806203}.

\bibitem[Lim and Tao(2025)]{LimTao2025}
K.~Lim and M.~Tao.
\newblock Appropriate state-dependent friction coefficient accelerates kinetic
  {L}angevin dynamics.
\newblock \emph{SIAM J. Appl. Math.}, 85\penalty0 (1):\penalty0 1--26, 2025.
\newblock \doi{10.1137/23M1625810}.

\bibitem[Lu(2026)]{lu2026sharp}
J.~Lu.
\newblock A sharp hypocoercive entropy decay estimate for underdamped
  {L}angevin dynamics.
\newblock \emph{arXiv preprint arXiv:2605.01933}, 2026.

\bibitem[Matthews et~al.(2018)Matthews, Weare, and Leimkuhler]{LMW}
C.~Matthews, J.~Weare, and B.~Leimkuhler.
\newblock Ensemble preconditioning for {M}arkov chain {M}onte {C}arlo
  simulation.
\newblock \emph{Stat. Comput.}, 28\penalty0 (2):\penalty0 277--290, 2018.
\newblock \doi{10.1007/s11222-017-9730-1}.

\bibitem[M{\"u}ller and Brown(1979)]{MuellerBrown1979}
K.~M{\"u}ller and L.~D. Brown.
\newblock Location of saddle points and minimum energy paths by a constrained
  simplex optimization procedure.
\newblock \emph{Theoret. Chim. Acta}, 53\penalty0 (1):\penalty0 75--93, 1979.
\newblock \doi{10.1007/BF00547608}.

\bibitem[Pavliotis(2014)]{pavliotis2014}
G.~A. Pavliotis.
\newblock \emph{Stochastic Processes and Applications: {D}iffusion Processes,
  the {F}okker-{P}lanck and {L}angevin Equations}, volume~60 of \emph{Texts in
  Applied Mathematics}.
\newblock Springer, New York, 2014.
\newblock \doi{10.1007/978-1-4939-1323-7}.

\bibitem[Roy et~al.(2018)Roy, Tan, and Flegal]{RoyTanFlegal2018}
V.~Roy, A.~Tan, and J.~M. Flegal.
\newblock Estimating standard errors for importance sampling estimators with
  multiple {M}arkov chains.
\newblock \emph{Statistica Sinica}, 28\penalty0 (2):\penalty0 1079--1101, 2018.
\newblock \doi{10.5705/ss.202016.0003}.

\bibitem[Sachs et~al.(2017)Sachs, Leimkuhler, and Danos]{SachsLeimkuhlerDanos}
M.~Sachs, B.~Leimkuhler, and V.~Danos.
\newblock {L}angevin dynamics with variable coefficients and nonconservative
  forces: {F}rom stationary states to numerical methods.
\newblock \emph{Entropy}, 19\penalty0 (12):\penalty0 647, 2017.
\newblock \doi{10.3390/e19120647}.

\bibitem[Sokal(1997)]{SokalMCLecture}
A.~D. Sokal.
\newblock {M}onte {C}arlo methods in statistical mechanics: {F}oundations and
  new algorithms.
\newblock In \emph{Functional Integration: Basics and Applications}, pages
  131--192. Springer, Boston, MA, 1997.
\newblock \doi{10.1007/978-1-4899-0319-8_6}.

\bibitem[Sundman(1913)]{sundman1913memoire}
K.~F. Sundman.
\newblock Memoir on the three-body problem.
\newblock \emph{Acta Mathematica}, 36\penalty0 (1):\penalty0 105--179, 1913.

\bibitem[Vilmart(2014)]{vilmart2014weak}
G.~Vilmart.
\newblock Weak second order multirevolution composition methods for highly
  oscillatory stochastic differential equations with additive or multiplicative
  noise.
\newblock \emph{SIAM Journal on Scientific Computing}, 36\penalty0
  (4):\penalty0 A1770--A1796, 2014.
\newblock \doi{10.1137/130935331}.

\end{thebibliography}
\end{document}